\documentclass[journal]{IEEEtran}
\IEEEoverridecommandlockouts
\usepackage{cite}
\usepackage{bm}
\usepackage{amsmath,amssymb,amsfonts}
\usepackage{algorithmic}
\usepackage{graphicx}
\usepackage{textcomp}
\usepackage{cleveref}
\usepackage{amsthm}
\usepackage{mathrsfs}
\usepackage{xcolor}
\def\BibTeX{{\rm B\kern-.05em{\sc i\kern-.025em b}\kern-.08em
    T\kern-.1667em\lower.7ex\hbox{E}\kern-.125emX}}

\usepackage{url}

\newtheorem{theorem}{Theorem}[section]
\newtheorem{proposition}{Proposition}[section]

\newtheorem{remark}{Remark}[section]
\newtheorem{lemma}{Lemma}[section]

\usepackage{mathtools}

\DeclarePairedDelimiter\floor{\lfloor}{\rfloor}

\newcommand{\blue}[1]{#1} 
\newcommand{\red}[1]{#1}
\newcommand{\orange}[1]{#1}

\DeclareMathOperator*{\argmin}{arg\,min}
    
\begin{document}
\title{A First-Order Gradient Approach for the Connectivity Optimization of Markov Chains}
\author{Christian~P.C.~Franssen, Alessandro~Zocca, and Bernd~F.~Heidergott
\thanks{C.P.C.~Franssen and B.F.~Heidergott are with the Department of Operations Analytics, VU Amsterdam, 1081 HV, Amsterdam, The Netherlands (email: c.p.c.franssen@vu.nl, b.f.heidergott@vu.nl).}
\thanks{A.~Zocca is with the Department of Mathematics, VU Amsterdam, 1081 HV, Amsterdam, The Netherlands (email: a.zocca@vu.nl).}}

\maketitle

\begin{abstract}
\red{Graphs} are commonly used to model various complex systems, including social networks, power grids, transportation networks, and biological systems. In many applications, the connectivity of these networks can be expressed through the Mean First Passage Times (MFPTs) of a Markov chain modeling a random walker on the graph.
In this paper, we generalize the network metrics based on Markov chains' MFPTs and extend them to networks affected by uncertainty, in which edges may fail and hence not be present according to a pre-determined stochastic model. To find optimally connected \blue{Markov chains}, we present a parameterization-free method for optimizing the MFPTs of the Markov chain. More specifically, \blue{we present an efficient Simultaneous Perturbation Stochastic Approximation (SPSA) algorithm} in the context of Markov chain optimization. The proposed algorithm is suitable for both fixed and random networks. Using various numerical experiments, we demonstrate scalability compared to established benchmarks. Importantly, our algorithm finds an optimal solution without requiring prior knowledge of edge failure probabilities, allowing for an online optimization approach.
\end{abstract}

\begin{IEEEkeywords}
Network analysis and Control, Markov processes, Optimization, Simultaneous Perturbation Stochastic Approximation, Optimization algorithms.
\end{IEEEkeywords}

\section{Introduction}
\label{sec:introduction}
\IEEEPARstart{G}{raphs} \red{are widely used to describe complex systems, such as social networks, power grids, transportation networks, and biological systems. }
In many applications, it is often useful to study the dynamical properties of these networked systems by analyzing the behavior of random walks on them. 
Such a random walk can be specified as a Markov chain by explicitly giving all transition probabilities between each pair of nodes.
Mean first passage times (MFPTs) characterize the connectivity properties of the given network, as they describe the expected time for the random walker starting at a given node to reach another specific node for the first time, for a review, see~\cite{Hunter2016,Hunter2018}. MFPTs are studied in various fields where phenomena can be modeled with Markovian dynamics, for instance, in biology~\cite{Chou2014} or chemistry~\cite{Kalantar2018,Weiss2007}.

In this paper, we focus on analyzing and minimizing mean first passage times (MFPTs) in Markov chains. \red{MFPTs serve as building blocks for several important metrics of network connectivity and robustness. Examples include the \textit{effective graph resistance} or \textit{Kirchhoff index}~\cite{klein1993resistance,ellens2011effective}, the \textit{DW-Kirchhoff index}~\cite{bianchi2019kirchhoffian}, and the \textit{Kemeny constant}~\cite{kemeny1976finite}, all of which can be expressed as weighted combinations of MFPTs.}

\red{A central theme in the literature is how to optimally design or modify edge weights to improve such metrics; see, e.g.,~\cite{ghosh2008minimizing,patel2015robotic}. In this context, optimizing edge weights can be equivalently interpreted as designing an optimal exploration policy over the graph. Indeed, since edge weights define the transition probabilities of the associated random walk, they implicitly determine how the graph is traversed. From this perspective, the problem becomes one of policy optimization, where the goal is to identify the Markov chain that minimizes a specified MFPT-based performance criterion.} 

The resulting problem is generally non-convex, but restricting the set of solutions to reversible Markov chains can turn the problem into a convex one. 
For example, the minimization of the Kemeny constant can be made convex by constraining the stationary distribution of the Markov chain and consequently can be solved using semi-definite programming \cite{patel2015robotic,Duan_Bullo_2020}.
However, the restriction to reversibility has a severe impact on the achievable objective value (see \cite{Duan_Bullo_2020} for an example), which motivates our research into a gradient-based optimization method that does not rely on convexity and thus allows one to drop the reversibility condition on the network.

In the present paper, we deviate from the standard approaches in the literature by formulating and solving more general and possibly non-convex \red{optimization problems for network connectivity}. 
In particular, we allow for more general MFPT-based objective functions.
Moreover, we extend the existing literature by considering settings in which network edges can fail with some (possibly correlated) probability.
For instance, this can be of particular interest when modeling power grids where transmission lines may fail due to climatic disasters (wildfires, floods, storms, etc.).
Therefore, in the rest of the paper, we distinguish between the \textit{fixed support case}, in which the set of edges of the network is fixed, and the \textit{random support case}, \blue{where there is uncertainty in the support of the graph, i.e., edges fail stochastically.}

Our main contributions are as follows:
\begin{itemize}
    \item We generalize the existing metrics for network connectivity based on MFPTs, accommodating both linear combinations and general (smooth) functions. \blue{Moreover, we introduce a broader class of network optimization problems with the new generalized connectivity metric as an objective function. We consider both the fixed support case and the non-trivial random support case.}

    \item If the considered graph is Hamiltonian, we derive explicit bounds for the sum of the MFPTs if we assume reversibility \blue{in Proposition~\ref{prop:ul_bound}}, \red{and we show that we can obtain a strictly better solution by relaxing the reversibility assumption.}

    \item We modify the Simultaneous Perturbation Stochastic Approximation (SPSA) algorithm \red{\cite{spall1992multivariate, Sadegh_1997} specifically for} the context of Markov chain optimization to solve possibly non-convex policy optimization problems, both in the fixed support setting \blue{(see Theorem~\ref{thm:est_sfd_conv})}, as well as in the random support setting \blue{(see Theorem~\ref{thm:est_sfd_ran})}.
    \red{A key feature of our approach is a tailored perturbation and projection scheme that guarantees all iterates remain valid Markov chains, thereby respecting the hard constraints intrinsic to the problem of Markov chain optimization.}
    Concurrently, we allow for constraining the stationary distribution of the Markov chain.
    We show that our method outperforms a benchmark in terms of scalability and can seamlessly integrate an online optimization approach.
    Our relatively simple yet powerful extension of SPSA underscores the versatility of first-order methods in studying network connectivity. 

    \item We demonstrate our approach in the context of surveillance applications taken from \cite{patel2015robotic}. Specifically, by relaxing the reversibility assumption, we can design much more efficient stochastic surveillance algorithms for the fastest detection of intruders and anomalies in networks, both in the fixed support setting and in the random support setting.
\end{itemize}

The paper is organized as follows. Section~\ref{sec:not} introduces the notation and some preliminaries.
In Section~\ref{sec:probdesc}, we describe the \red{Markov chain} optimization problem with the generalized MFPT metric as objective in the fixed support setting and provide conditions for its convexity. In the same section, for some special graphs of interest, we also derive explicit bounds for the so-called ``price of reversibility''. 
We derive closed-form expressions for the derivative of MFPTs with respect to the parameters of the Markov chain in Section~\ref{sec:dif}, and leverage them to introduce an SPSA-type method for the optimization of Markov chains. 
In Section~\ref{sec:ran}, we consider the problem in the random support setting, explicitly list sufficient conditions for convexity, and explain how to modify our algorithm accordingly.
In Section~\ref{sec:ne}, we demonstrate SPSA's scalability and its handling of correlated edge presence in the random support setting.
We present an application of SPSA in a surveillance context in Section~\ref{sec:ne_2}.
Using our SPSA, we obtain a solution that, leveraging directionality and irreversibility, strictly improves reversible solutions. We conclude with a discussion of further research in Section~\ref{sec:conclusion}.

\section{Notation and preliminaries} \label{sec:not}
We consider a finite, directed, and simple graph $ \red{\mathscr{G}}  = (V, E) $, with set of nodes 
$ V= \{ 1 , \ldots , N \}$ and set of directed edges $ E \subset V \times V$. We denote by $E_i$ the set of edges leaving node $i$, i.e., $E_i = \{\ell \in E : \ell=(i,j ), j \in V \}$.
Given $\red{\mathscr{G}}=(V,E)$, we consider the vector $x \in [0,\infty)^{|E|}$ whose elements are the edge weights of $\red{\mathscr{G}}$, that is, $x_{\ell} = x_{(i,j)}$ is the weight of edge $ \ell = (i,j) \in E$.
We denote the vector of weights of the edges leaving node $i \in V$ by
\begin{align} \label{eq:xi}
    x_{(i,\cdot)} = \left(x_{(i,j)}\right)_{(i,j) \in E_i},
\end{align}
and, consequently, the vector $x$ can be seen as the concatenation of these vectors, i.e.,
\begin{align} \label{eq:vec}
    x = \begin{bmatrix} x_{(1,\cdot)} & \dots & x_{(N,\cdot)} \end{bmatrix}.
\end{align}

Throughout the paper, we specifically consider weights that induce a Markov chain on a graph $\red{\mathscr{G}}$.
More formally, let $ \mathcal{P} = \mathcal{P}(\red{\mathscr{G}})$ be the set of stochastic matrices associated with $\red{\mathscr{G}}$, i.e., the collection of \red{$N \times N$} matrices $P$ such that $P_{ij} \geq 0$ if $(i,j) \in E$ and $P_{ij} = 0$ otherwise. 
Then, for $x\in [0,\infty)^{|E|}$ and $\sum_{\ell \in E_i} x_\ell > 0$, for all $i \in V$, we define the Markov chain associated with $x$ through  
\begin{align} \label{eq:x_to_P}
    P(x) = P(x, \red{\mathscr{G}}) = \sum_{\ell=(i,j) \in E} \frac{x_{\ell}}{\sum_{\ell \in E_i} x_{\ell}} e_i e_j^{\top} \in \mathcal{P},
\end{align}
where $e_i \in \mathbb{R}^{N}$ is the vector with a 1 in the $i$th position and zeros elsewhere.
For our analysis, we work with the following constraint sets:
\begin{align*}
\begin{array}{rrl}
     & \mathcal{X}^{\mathrm{pos}} &= \Big\{ x \in \mathbb{R}^{|E|} : 
x_\ell \geq 0 , \ \forall \ \ell \in E \Big\} , \\
     & \mathcal{X}^{(\mathrm{pos},\varepsilon)} &= \Big\{ x \in \mathbb{R}^{|E|} : 
x_\ell \geq \varepsilon , \ \forall \ \ell \in E \Big\} ,
\end{array}
\end{align*}
for $ \varepsilon > 0 $ small, and
\[
\mathcal{X}^{\mathrm{eq}} = \Big\{ x \in \mathbb{R}^{|E|} : \sum_{\ell \in E_i} x_{\ell} = 1, \ \forall \ i \in V \Big\} .
\]
Let 
\[
\mathcal{X} = \mathcal{X}^{\mathrm{pos}} \cap \mathcal{X}^{\mathrm{eq}}
\text{, and } \mathcal{X}_\varepsilon = \mathcal{X}^{(\mathrm{pos},\varepsilon) } \cap \mathcal{X}^{\mathrm{eq}}.
\]
With this notation, (\ref{eq:x_to_P}) is a Markov chain for $ x \in \mathcal{X}^{\mathrm{pos}}$, if $\sum_{\ell \in E_i} x_\ell > 0$, for all $i \in V$. Note that the latter is implied for $ x \in {\cal X}^{\mathrm{eq}}$.
Therefore, it is natural to encode the normalization in \eqref{eq:x_to_P} into the set of feasible weights via $ \mathcal{X}^{\mathrm{pos}}$.
Then, for $x \in \mathcal{X}$, \eqref{eq:x_to_P} simplifies to
\begin{align} \label{eq:x_to_P2}
    P(x) = P(x, \red{\mathscr{G}}) = \sum_{\ell=(i,j) \in E} x_{\ell} e_i e_j^{\top} \in \mathcal{P}.
\end{align}
Throughout the paper, we will refer to a general stochastic matrix associated with $\red{\mathscr{G}}$ as $P$ if there is no explicit dependence on $x$.

The Markov chain associated with a given stochastic matrix $P \in \mathcal{P}(\red{\mathscr{G}})$ describes a random walk on the graph $\red{\mathscr{G}}$. More formally, if $X_n \in V$ denotes the node of the network at which the random walker resides at time $n \in \mathbb{N}$, we assume that for every $i,j \in V$, the transition probabilities from node $i$ to node $j$ only depend on the current node and are given by the corresponding element of the matrix $ P $, i.e.,
\begin{align*}
\begin{split}
    P_{ij} &= \text{Prob}(X_{n+1}=j ~|~ X_{n}=i )  \\
    &= \text{Prob}(X_{n+1}=j ~|~ X_{n}=i, \dots, X_{0}).
\end{split}
\end{align*}
Then, $ \{ X_n \}_{n \in \mathbb{N}} $ is a (homogeneous) Markov chain with node transition matrix $ P $.

Throughout the paper, we assume that $ \red{\mathscr{G}} $ is strongly connected. This assumption implies that any $ P (x)  $ is irreducible whenever $ x \in \mathcal{X}_\varepsilon $. \blue{In an irreducible Markov chain, all nodes in the graph ``communicate'', i.e., for all $(i,j) \in V\times V$, there exists an index $k \in \mathbb{N}$, such that $P_{ij}^k > 0$.}
\orange{An irreducible and finite $ P(x) $ is positive recurrent and so has a unique stationary distribution $ \pi= \pi ( P(x) )  $, such that 
\begin{align*}
    \pi > 0, \ \sum_{i=1}^N \pi_i = 1, \text{ and } \pi P(x) = \pi.
\end{align*} 
Throughout the paper, we adopt the convention that $ \pi(P(x)) $ (or simply $\pi$ when $P$ is clear from context) is a row vector and always refers to this unique stationary distribution. Moreover, we denote
\begin{align*}
    \blue{\mathcal{X}^{\hat{\pi}} = \{ x \in \mathcal{X} : \pi(P(x)) = \hat{\pi} \}, \text{ and } \red{\mathcal{X}^{\hat{\pi}}_\varepsilon} =  \mathcal{X}^{\hat{\pi}} \cap \mathcal{X}_\varepsilon,}
\end{align*}
as sets inducing Markov chains} \blue{$ \red{\mathcal{P}^{\hat{\pi}}} = \{ P \in \mathcal{P} : \pi(P) = \hat{\pi}\}$} with \blue{fixed} stationary distribution $\hat{\pi}$ via \eqref{eq:x_to_P2}.

The value \red{$ \pi_j$} of the stationary distribution of $ P $ can be interpreted as the long-run fraction of visits to node $ \red{j}$ by the random walker.
\orange{A convenient way to formalize this is via the ergodic projector of $P$, that is, the matrix $ \Pi := ( \red{\Pi_{i  j }} )_{i,j \in V},$
whose entries are defined as
\begin{eqnarray*}
\red{\Pi_{i  j }} := \lim_{n \to \infty} \frac{1}{n} \sum_{m=1}^n \text{Prob}(X_{m}=j ~|~ X_{0}=i).
\end{eqnarray*}
In particular, if $P$ is irreducible (but not necessarily aperiodic), for any $i \in V$ it holds a.s.~that
\begin{eqnarray*}
\Pi_{i  j }  = \lim_{n \to \infty} \frac{1}{n} \sum_{m=1}^n {\bf 1 }( X_m = \red{j} ) ,
\end{eqnarray*}
where ${\bf 1 }(\cdot) $ denotes the indicator function.}
\red{In other words, the limits $ \Pi_{i  j } $ }are independent of the initial node $ i $ (i.e., $ \red{\pi_{j}} = \red{\Pi_{i  j }}$, for all $ j \in V$), so that $ \Pi $ has equal rows, all equal to the unique stationary distribution \red{$ \pi = ( \pi_j : j \in V )$} of $ P $.

For a Markov chain $P$, let $ M \in \red{\mathbb{R}_{\geq 0}^{N\times N}}$ be the matrix of the mean first passage times, where for every $i,j \in V$ $$M_{ij} = \mathbb{E}\left[ \min\{ n \geq 1 ~|~ X_n = j, X_0 = i \}\right]$$ denotes the mean number of steps it takes a random walker to go from $i$ to $j$. 
From the irreducibility of $P$ and the finiteness of the set of nodes $V$, it follows that $ M_{ i j } < \infty $ for every $ i, j  \in V$. 
The entries of $ M $ can be interpreted as ``distances'' or ``connectivity metrics'' between nodes: if $ M_{ i j } < M_{ i k }$, then node $ j $ is reachable from node $ i$ on average in fewer steps than node $ k$.
The {\em deviation matrix} of $P$ is defined as 
\begin{equation}\label{eq:matrixDefDevMatrix}
D = (I - P + \Pi)^{-1} - \Pi.
\end{equation}
Given $ \Pi $ and $ D $, the mean first passage time matrix $M$ can be obtained using the following closed-form expression \cite{kemeny1976finite, berkhout2019analysis}:
\begin{equation}\label{eq:m1}
M = (I - D + {\bar 1} {\bar 1}^\top \textup{dg}(D)) \cdot \textup{dg}(\Pi )^{-1},
\end{equation} 
\blue{where $\bar{1}$ is an appropriately sized vector of ones and $\textup{dg}(\cdot)$ denotes the diagonal matrix formed by placing the elements of the argument matrix on the diagonal.}

A Markov chain $P$ is \textit{reversible} if it satisfies the following condition
\begin{itemize}
\item [\textbf{(R)}] 
{\em $P$ satisfies the detailed balance equations 
\begin{align*}
\pi_i(P) P_{ i j } = \pi_j(P) P_{ j i}, \quad \forall \,i , j \in V,
\end{align*}}
\end{itemize}
see \cite{aldous2014reversible}.
\blue{In the following, we let $\mathcal{P}^{\mathrm{rev}} = \{ 
P \in \mathcal{P} :  \pi_i(P) P_{ i j } = \pi_j(P) P_{ j i}, \ \forall \,i , j \in V\}$ denote the set of reversible Markov chains.
In Lemma~\ref{lem:surj}, we will show that $P(x)$ via \eqref{eq:x_to_P} is a surjective map from 
\begin{align*}
\mathcal{X}^{(\mathrm{sym}, 1)} = \mathcal{X}^{\mathrm{sym}} \cap \left\{ x \in \mathcal{X}^{\mathrm{pos}} : \sum_{(i,j) \in E}x_{(i,j)} = 1 \right\},
\end{align*}
where $\mathcal{X}^{\mathrm{sym}} = \left\{ x \in \mathcal{X}^{\mathrm{pos}} : x_{(i,j)} = x_{(j,i)}, \, \forall \ (i,j)  \in E\right\}$, to $\mathcal{P}^{\mathrm{rev}}$. Furthermore, $P(x)$ becomes a bijection if we only consider strongly connected graphs.}

\begin{lemma} \label{lem:surj}
    \blue{{\em $P(x)$ via \eqref{eq:x_to_P} is a surjective mapping from $\mathcal{X}^{(\mathrm{sym}, 1)}$ to $\mathcal{P}^{\mathrm{rev}}$. Furthermore, $P(x)$ is a bijection when $\mathcal{X}^{(\mathrm{sym}, 1)}$ is restricted to induce strongly connected graphs.}}
    
    \begin{proof}
        \blue{See Appendix~\ref{sec:proof_surj}.}
    \end{proof}
\end{lemma}

A frequently used quantity \red{to study network connectivity} is the \textit{effective graph resistance} or \textit{Kirchhoff index}.
For a graph $(\red{\mathscr{G}}, x)$, $x\in \mathcal{X}^{\mathrm{sym}}$, the effective graph resistance $R_{tot}(x)$ is defined using the effective resistance matrix $R(x) = R(\red{\mathscr{G}},x) \in \mathbb{R}_{\geq 0}^{N\times N}$ as follows 
\begin{align} \label{eq:egr_pw}
 R_{tot}(x) = R_{tot}( \red{\mathscr{G}}, x ) = \sum_{i < j} R_{ij}(x),
\end{align}
where the effective resistance $R_{ij}(x)$ between node $i$ and $j$ can be expressed in terms of the mean first passage matrix $M$ corresponding to the Markov chain $P(x)$ as
\begin{align} \label{eq:egr}
R_{ij} (x) = \frac{1}{\sum_{\ell \in E} x_{\ell}} \left( M_{ i j } + M_{ji} \right),
\end{align}
where $ R(x)$ is symmetric by construction \cite{ghosh2008minimizing}.
If the effective graph resistance is small, then a random walker on the corresponding graph quickly reaches any node from any other node.
Therefore, $ R_{tot}(x) $ is a measure of connectivity, with lower values of $ R_{tot} (x)$ corresponding to higher connectivity.
The effective graph resistance is used as a proxy for robustness in the design and control of power transmission networks (see~\cite{Zocca2021}): a low effective graph resistance indicates that all nodes in the network are well connected so that, in the event of transmission line failures, electricity can still efficiently flow throughout the surviving network \cite{wang2015network,dorfler2010synchronization}.
This concept extends naturally to communication networks~\cite{tizghadam2010betweenness, rueda2017robustness}, and air transportation networks~\cite{yang2018designing}.

\section{Markov chain connectivity on fixed graphs}
\label{sec:probdesc}
Let us consider a generic metric of the connectivity of a Markov chain $P$ expressed as the (weighted) sum of MFPT's, defined as 
\begin{align} \label{eq:obj}
 S (P , C ) =\sum_{i,j} C_{ i j }(P) M_{ i j },
\end{align}
where $ C(P) \in \mathbb{R}_{\geq 0}^{N \times N}$ is a matrix of non-negative weights $C_{ij}(P) \geq 0$ that may depend on the Markov chain $P$ itself.
One may choose these weights such that $ C_{ i j }(P) >  C_{ k l }(P) $ to express the fact that the decision-maker values the connectivity of $ i $ to $j $ more than the connectivity from $ k  $ to $ l $. 
\blue{In this way, it becomes possible to focus the network optimization problem around a subset of nodes or specific clusters of the given network.}
A standard choice for $C(P)$ is 
\begin{align} \label{eq:c_kem}
C_{ij}^{\pi}(P) = \pi_i(P) \pi_j(P),
\end{align}
which yields the \textit{Kemeny constant} $ K = S (P , C^\pi ) $. The Kemeny constant is a measure of the strength of connectivity of a Markov chain and is frequently used in graph analysis (e.g., see \cite{berkhout2019analysis,yilmaz2020kemeny,patel2015robotic}).

Another canonical choice is to take $C(P)$ as the scalar matrix $C^{\bar{1}} = \bar{1}\bar{1}^{\top} - I$. If $\sum_{\ell \in E} x_{\ell} = 1$, we recover the \textit{DW-Kirchhoff} (directed \& weighted) index introduced by~\cite{bianchi2019kirchhoffian} as a generalization of the effective graph resistance \eqref{eq:egr_pw}.

In this paper, our goal is to find optimally connected Markov chains by minimizing $S (P , C )$. From now on, we will treat $C$ as a scalar matrix; however, it is important to note that our optimization techniques also apply to cases where $C$ varies with $P$ \red{(as long as $C(P)$ is three times differentiable)}, for example, when optimizing the graph's Kemeny constant.

\subsection{Optimizing Markov chain connectivity}

Let us formulate the problem of designing an optimally connected Markov chain $P \in \mathcal{P}$ on a graph $\red{\mathscr{G}}$ as
\begin{align}\label{eq:main_problem}
\blue{\min_{P \in {\cal P} }}\quad &  \blue{S ( P , C ).}
\end{align}
We also consider the variant of the optimization problem~\eqref{eq:main_problem} with an additional constraint specifying a target stationary distribution $ \hat \pi \in (0,1)^N$ for $P$, that is, 
\begin{align}\label{eq:main_problem2}
    \begin{array}{rl}
    \blue{\displaystyle \min_{P \in {\cal P}^{\hat{\pi}} }} &  \blue{S ( P , C ).}
    \end{array}
\end{align} 
This version of the problem is studied in \cite{patel2015robotic} in the context of finding a stochastic surveillance policy $P$ for network surveillance, where the policy is designed to visit nodes in proportion to the anticipated probability of intruder appearances at each node, represented by $\hat{\pi}$.

If $P$ is reversible and $C = C^{\hat{\pi}} = ( \Hat{\pi}_i \Hat{\pi}_j )_{i,j \in V}$, for fixed stationary distribution $\Hat{\pi}$, problem \eqref{eq:main_problem2} has been shown to be convex \cite{patel2015robotic}.
\blue{Note that dropping the $\pi(P) = \hat{\pi}$ constraint causes the problem to become non-convex, since the set of reversible Markov chains is not convex due to the presence of non-linearity in condition (\textbf{R}).
However, in the following theorem, we show how optimization via $\mathcal{X}^{(\mathrm{sym},1)}$ preserves convexity when generalizing 
\eqref{eq:main_problem2} to \eqref{eq:main_problem}, thus relaxing the $\pi(P) = \hat{\pi}$ constraint while also allowing for any non-negative and symmetric $C$.}

\begin{theorem} \label{thm:convexity}
    \blue{{\em Let $C \in \mathbb{R}_{\geq 0}^{N \times N}$ be symmetric. Consider the optimization problem
    \begin{equation} \label{eq:main_problem_rev}
    \begin{aligned}
        \min_{x \in {\cal X}^{(\mathrm{sym},1)} }\quad &  S ( P(x) , C ). \\
    \end{aligned}
    \end{equation}
    Then, the reversibility condition \textbf{(R)} is satisfied for $P(x)$ via \eqref{eq:x_to_P} for all strongly connected $(\red{\mathscr{G}},x)$ 
    and the problem in~\eqref{eq:main_problem_rev} is convex.
    Moreover, if $C$ is a matrix with identical entries $C_{ij} = c > 0$, then \eqref{eq:main_problem_rev} is a strictly convex problem.}
    }
    
    \begin{proof}
        \blue{As argued earlier, we can model all $P$ that satisfy condition \textbf{(R)} 
        using $\mathcal{X}^{(\mathrm{sym}, 1)}$, which is clearly convex in view of its linear constraints.
        To show the convexity of $S ( P(x) , C )$, we use the convexity of all pairwise resistance distances $R_{ij}(x)$ in any undirected weighted graph $(\red{\mathscr{G}},x)$ shown in \cite{ghosh2008minimizing}.       
        It follows that for any symmetric matrix $C \in \mathbb{R}_{\geq 0}^{N \times N}$ we can rewrite
        \begin{align*}
        S(P(x),C) &= \sum_{i < j} C_{ij} ( M_{ij} + M_{ji} ) =  \sum_{i < j} C_{ij} R_{ij}(x).
        \end{align*}
        Thus, being the sum of convex functions, $S(P(x),C)$ is convex.}

        \blue{To show strict convexity, assuming $C$ with identical entries $c > 0$, we use the fact that the effective graph resistance $R_{tot}(x) = \sum_{i < j} R_{ij}(x)$ is strictly convex \cite{ghosh2008minimizing}. Indeed, we can simply rewrite
        \begin{align*}
        S(P(x),C) =  \sum_{i<j} c R_{ij} (x) = c\sum_{i<j} R_{ij}(x),
        \end{align*}
        and strict convexity follows when $c > 0$.}
        \end{proof}
\end{theorem}

\begin{remark} \label{rem:infty}
    \blue{We do not need to explicitly model the constraint that $(\red{\mathscr{G}},x)$ is strongly connected in Theorem~\ref{thm:convexity}, since any graph, where $C_{ij} > 0$ and $i$ and $j$ are in separate components, has $S ( P(x) , C ) = \infty$ (due to $M_{ij} = \infty$).}
\end{remark}

\blue{Note that in the more general case where the matrix $C$ is symmetric but contains non-identical entries, problem \eqref{eq:main_problem_rev} is only convex but not necessarily strictly convex.
Moreover, the findings in Theorem~\ref{thm:convexity} carry over to problem \eqref{eq:main_problem2} since enforcing a fixed stationary distribution $\pi(P) = \hat{\pi}$ eliminates the non-linearity in (\textbf{R}), and so $\mathcal{P}^{\hat{\pi}}$ is convex.}  

\blue{Condition \textbf{(R)} makes the optimization problem~\eqref{eq:main_problem} convex, but adding it also dramatically restricts the set of possible solutions. 
In the next subsection, we show the negative impact of the reversibility condition on the optimal solution for $S(P,C)$, providing a strong rationale for addressing the problem \textit{without} assuming reversibility:}
\begin{align}\label{eq:main_problem_non_rev}
    \blue{\min_{x \in {\cal X} }} \quad &  \blue{S ( P(x) , C ),}
\end{align}
where $P(x)$ via \eqref{eq:x_to_P2}, and the analogous problem for \eqref{eq:main_problem2}:
\begin{align}\label{eq:main_problem2_non_rev}
    \blue{ \min_{x \in \mathcal{X}^{\hat{\pi} }} \quad} & \blue{ S ( P(x) , C ).}
\end{align}
\blue{Note that these problem classes can be extended to cases where $C$ depends on $P$, as is the case, for instance, in the minimization of the Kemeny constant, see~\eqref{eq:c_kem}.}

\subsection{The price of reversibility} \label{sec:por}

\blue{In the following, we consider a specific case that allows for comparing analytical reversible and non-reversible solutions to \eqref{eq:main_problem}.
To that end, consider a directed graph $\red{\mathscr{G}}$ that is \textit{Hamiltonian}, which means that $\red{\mathscr{G}}$ admits a path, called \textit{Hamiltonian cycle}, that starts and ends at the same node, visiting all other nodes in the graph exactly once. 
The following theorem states that if $\red{\mathscr{G}}$ is a Hamiltonian graph and if we assume unit weights $ C = C^{\bar{1}}$ for the mean first passage times in \eqref{eq:obj}, then the optimal solution of~\eqref{eq:main_problem} yields a Hamiltonian cycle on $ \red{\mathscr{G}} $.}

\begin{proposition}\label{prop:tour}{\em 
Let $ \red{\mathscr{G}} = ( V, E )$ be a directed Hamiltonian graph. Then, the Markov chain 
\begin{equation}  \label{eq:opt_P_nonrev}
    \blue{\Hat{P} = \argmin_{P \in \mathcal{P}} \ S ( P ,  C^{\bar{1}}  ),}
\end{equation}
corresponds to a Hamiltonian cycle on $\red{\mathscr{G}}$, that is, $ \Hat{P} $ is irreducible, and each row has one entry equal to $ 1 $ and all other entries equal to $ 0 $.
Furthermore, its optimal value equals
\begin{align}
    S ( \Hat{P},  C^{\bar{1}}  ) = \frac{N^3 - N^2}{2}.
\end{align}
}   

\begin{proof}
The proof straightforwardly follows from \cite[Theorem 2]{Borkar_Miclo_2020}, \blue{see Appendix~\ref{sec:proof_tour}.}
\end{proof}
\end{proposition}

If the graph $\red{\mathscr{G}}$ admits a Hamiltonian cycle and we assume reversibility, the following lower and upper bounds hold for the optimal Markov chain $\hat{P}$.
\begin{proposition}\label{prop:ul_bound}
{\em Let $ \red{\mathscr{G}} = ( V, E )$ be a directed Hamiltonian graph with $N \geq 3$ nodes. Consider the optimization problem~\eqref{eq:opt_P_nonrev} with $P$ satisfying \textbf{(R)} as an additional constraint;
\begin{align}  \label{eq:opt_rev}
    \blue{\check{P} = \argmin_{P \in \mathcal{P}^{\mathrm{rev}}}} \ \blue{S ( P ,  C^{\bar{1}}  ). } 
\end{align}
Then, the following inequality holds for its optimal value:
\begin{align} \label{eq:lu_bound}
   N^3 - 2N^2 + N \leq S(\check{P}, C^{\bar{1}}) \leq \frac{N^4 - N^2}{6}.
\end{align}
}

\begin{proof} 
See Appendix~\ref{sec:proof_ub_lb}.
\end{proof}
\end{proposition}

Inspecting \cref{prop:tour,prop:ul_bound}, we conclude that one always obtains strictly worse solutions by assuming reversibility, since $\frac{N^3 - N^2}{2} < N^3 - 2N^2 + N$
for any $N \geq 3$.

\subsection{Illustration of the price of reversibility}
\begin{figure}[h!]
    \centerline{
    \includegraphics[width=0.585\columnwidth]{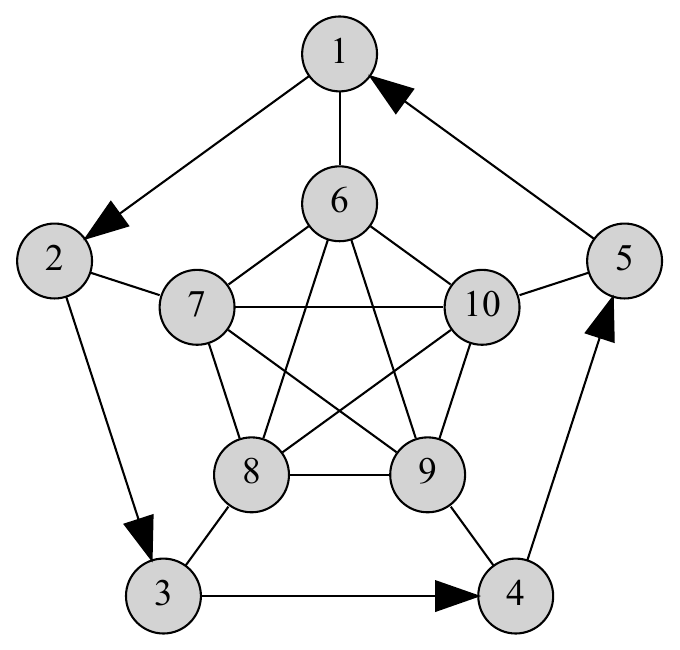}}
    \caption{Directed Petersen Graph, $N = 10$.}
    \label{fig:double_star}
\end{figure}

\begin{figure}[h!]
    \centerline{
    \includegraphics[width=0.585\columnwidth]{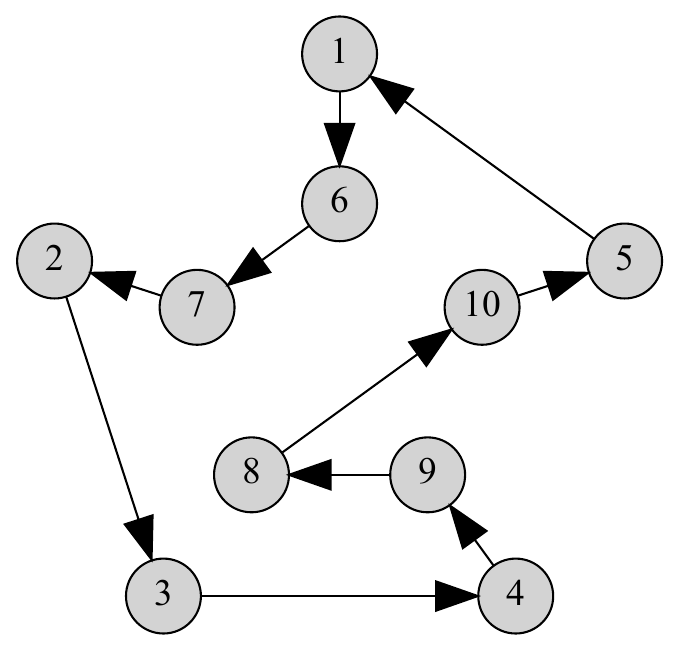}}
    \caption{Global optimal non-reversible solution $\hat{P}$, which yields an optimal value of $S(\hat{P}, C^{\bar{1}}) = 450$.}
    \label{fig:double_star_non_rev}
\end{figure}

We demonstrate the price of reversibility in the case of a directed version of the Petersen graph, see Fig.~\ref{fig:double_star}. 
Here, the edges with missing arrows can be traversed in both directions.
For this instance, the optimal non-reversible Markov chain, see \eqref{eq:opt_P_nonrev}, is visualized in Fig.~\ref{fig:double_star_non_rev}.

\begin{figure}[h!]
     \centerline{
    \includegraphics[width=0.585\columnwidth]{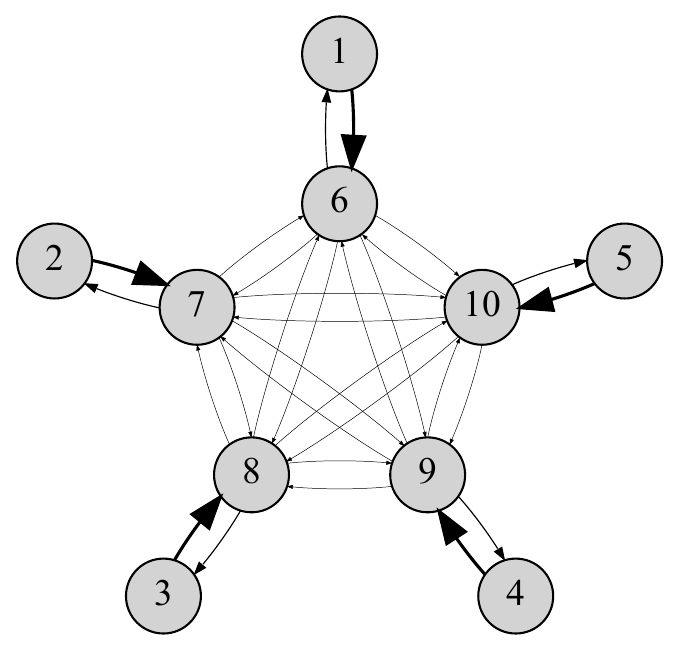}}
        \caption{Global optimal reversible solution $\check{P}$, which yields an optimal value of $S(\check{P}, C^{\bar{1}}) = 1529$.}
        \label{fig:double_star_rev}
\end{figure}

This graph allows for multiple tours, so the global minimum is not unique.
For sake of comparison, we show the optimal reversible Markov chain as in \eqref{eq:opt_rev}, after having excluded the one-directional edges of the graph, see Fig.~\ref{fig:double_star_rev}, since these edges cannot be used in a reversible $P$.

\section{First-order derivatives and approximations} \label{sec:dif}

Simultaneous Perturbation Stochastic Approximation (SPSA) is a gradient-based optimization technique that is particularly effective in scenarios involving complex, noisy, and high-dimensional problem spaces. 
The key ingredient of SPSA is that the gradient of the objective function is approximated by using only two measurements, regardless of the dimensionality of the parameter space. 
Thus, SPSA significantly reduces the computational burden compared to traditional finite-difference methods \cite{spall1992multivariate, shi2021sqp}.
Working towards our final SPSA-based algorithm, we first provide analytical expressions for the derivatives of Markov chains. Subsequently, we show how to extend SPSA to approximate these derivatives.

Consider the generic optimization problem $$ \min_{x \in \mathbb{R}^d} J(x),$$ where the function $J$ is assumed to be three times continuously differentiable with respect to $ x $.
The SPSA gradient proxy of $J$ in $x$ \cite{spall1992multivariate} is obtained as follows: 
\begin{align} \label{eq:spsa_est}
    G(x) = & \ G(x, \Delta, \eta) = \  \frac{J(x + \eta \Delta) - J(x- \eta \Delta)}{2\eta} \Delta, 
\end{align}
\red{where $\Delta \sim \text{Rad}^d$ is a perturbation vector with entries independently and uniformly distributed over $\{-1, 1\}$.}
We can find a stationary point of $J$ using the following basic stochastic approximation recursion initiated at some initial point $x(0)$, 
\begin{align} \label{eq:spsa_rec}
    x(k+1) = x(k) - \alpha(k) G_k(x(k)),
\end{align}
for $k\in \mathbb{N}$, with step size sequence $\{\alpha(k)\}_{k \in \mathbb{N}}$ and gradient proxy $G_k(x(k)) = G(x(k), \Delta(k), \eta(k))$ depending on sequences $\{\Delta(k)\}_{k \in \mathbb{N}}$, and $\{\eta(k)\}_{k \in \mathbb{N}}$.
The algorithm converges strongly to a point $x$ for which $\nabla J(x) = 0$ if appropriate sequences $\{\alpha(k)\}_{k \in \mathbb{N}}$, $\{\Delta(k)\}_{k \in \mathbb{N}}$, and $\{\eta(k)\}_{k \in \mathbb{N}}$ are used (for a convergence result, see \cite{spall1992multivariate}).

\red{As we show in Appendix~\ref{sec:inf_spsa}, standard SPSA is not applicable for solving problem \eqref{eq:main_problem_non_rev}, since the constraints imposing $\mathcal{X}$ are \textit{hard constraints} and therefore cannot be violated.
Therefore, we adapt the generic SPSA algorithm in (\ref{eq:spsa_rec}) to operate within the constraint set $ \mathcal{X}$.
\orange{To that end}, we derive in Section~\ref{sec:ander} an analytical expression for the steepest feasible descent on $\mathcal{X}^{\mathrm{eq}}$.
Then, in Section~\ref{sec:spsa1}, we introduce a version of SPSA that estimates this descent direction, ensuring that $ x \pm \eta \Delta  \in \mathcal{X}^{\mathrm{eq}}$.
\orange{Independently, Section~\ref{sec:spsa2} presents an SPSA variant following~\cite{Sadegh_1997} that is designed to maintain the feasibility of the perturbed values $x \pm \eta \Delta$ with respect to the inequality-constrained set $\mathcal{X}^{\mathrm{pos}}$ across all iterates of $x$. To ensure this, we require that each iterate $x$ is projected back onto the restricted set $\mathcal{X}^{(\mathrm{pos},\varepsilon)}$, so that the resulting perturbations remain within $\mathcal{X}^{\mathrm{pos}}$.
In Section~\ref{sec:spsa3}, we then combine these two components, with equality feasibility from Section~\ref{sec:spsa1} and inequality feasibility from Section~\ref{sec:spsa2}, into a SPSA algorithm that guarantees $x \pm \eta \Delta$ is feasible for any $x \in \mathcal{X}_\varepsilon$, thereby solving problem~\eqref{eq:main_problem_non_rev} over $\mathcal{X}_\varepsilon$.}  
Lastly, we discuss how to apply SPSA to the extended problem~\eqref{eq:main_problem2_non_rev} in Section~\ref{sec:spsa4}, where a similar approach is used to optimize over $\mathcal{X}^{\hat{\pi}}_\varepsilon$.}

\subsection{Analytical derivatives} \label{sec:ander}

\subsubsection{Directional derivatives}
Denote the set of all feasible directions leaving \red{$x \in \mathcal{X}_\varepsilon $ with respect to the set $\mathcal{X}$} as
\[
 \mathcal{D} ( x; \mathcal{X} ) = \{  \delta \in \mathbb{R}^{|E|} : \: \| \delta \|_2 \leq 1 \land \exists \eta :   x +\eta \delta \in {\cal X } \}.
\]
Then, we can define the directional derivative for a continuously differentiable function $J(x)$ in the feasible direction $\delta \in \mathcal{D}  ( x ; \mathcal{X})$ as
\begin{align*}
\partial_\delta J(x) = \lim_{\eta \to 0} \frac{J(x+ \eta \delta)- J(x)}{\eta} .
\end{align*}
In Theorem~\ref{th:grad}, we provide analytical directional derivatives for $\Pi$, $D$, and $M$.

\begin{theorem}\label{th:grad}
{\em Let \red{$x\in \mathcal{X}_\varepsilon$} such that $P(x) \in \mathcal{P}$ and take \red{$\delta \in \mathcal{D}(x, \mathcal{X})$}.
The directional derivative of $ \Pi $ with respect to $ \delta $ is given by
\[
\partial_{\delta} \Pi
 = \Pi P'  D ,
\]
and that of $ D $ by
\[
\partial_{\delta} D =
- \Pi D + D P' D,
\]
where we use the abbreviation $P'= \partial_\delta P(x)$.
Moreover, for $ M = M ( P(x) ) $ the following identity holds 
\begin{eqnarray*}
 & &  \! \! \! \!  \! \! \! \! \partial_{\delta} M
  =  \left( -D' + \bar{1}\bar{1}^{\top} \textup{dg}(D') \right)\cdot \textup{dg}(\Pi)^{-1} \\
    & & + \left( I - D + \bar{1}\bar{1}^{\top} \textup{dg}(D) \right) \cdot \left( - \textup{dg}(\Pi)^{-1} \textup{dg}(\Pi') \textup{dg}(\Pi)^{-1}\right) ,
\end{eqnarray*}
with $ \Pi'= \partial_{\delta} \Pi $ and $ D'=\partial_{\delta} D $.
}

\begin{proof}
The differentiation of $ \pi$ with respect to $ P$ reduces the problem to the well-studied problem of differentiation of parameterized Markov chains.
The explicit solution follows, e.g., from \cite{heidergott2003taylor}, where it is also shown that $ \pi $ is analytic in the entries of $P$.
For the differentiability of $D$, we refer to \cite{leder2010approximation}, while that of $M$ follows from (\ref{eq:m1}) using the fact that $(A^{-1})' = -A^{-1} A' A^{-1}$.
\end{proof}
\end{theorem} 

Following Theorem~\ref{th:grad}, we can explicitly solve for directional derivatives of $S(P(x),C)$ in direction $\delta \in \mathcal{D}  ( x ; \mathcal{X})$:
\[
\partial_\delta S(P(x),C) = \sum_{i,j} C_{ij} (\partial_{\delta} M)_{ij}.
\]

\subsubsection{Steepest descent direction}
We will now show how to analytically find the steepest feasible direction $\delta(x) \in \mathcal{D} \left(x; \mathcal{X}\right)$ provided $x \in \mathcal{X}_\varepsilon$ and for $\varepsilon > 0$ small, that is, for $J(x) = S(P(x),C)$ we solve
\begin{align} \label{eq:sfd}
    \delta(x) = \argmin_{\delta \in \mathcal{D} \left(x; \mathcal{X}\right)} \ \partial_\delta J(x) .
\end{align}
For any $x \in \mathcal{X}_\varepsilon$, $x$ is strictly feasible with respect to non-negativity constraints of $\mathcal{X}$, i.e., $x > 0$.
Therefore, a feasible direction $\delta \in \mathcal{D} (x; \mathcal{X})$ preserves that $P(x+\eta \delta) \in \mathcal{P}$, for $\eta$ sufficiently small.
We can derive the necessary conditions for $\delta$ by writing the equality constraints in $\mathcal{X}$ as a system of equations $Ax = b$, where
\begin{align} \label{eq:A}
A ^{N \times |E|}= \begin{bmatrix}
    \Bar{1}^{\top}_{|E_1|} & & \\
     & \ddots & \\
     & & \Bar{1}^{\top}_{|E_N|}
\end{bmatrix}, \text{ and }
 b = \Bar{1}.
\end{align}
It follows that for any feasible direction $\delta\in \mathcal{D}(x; \mathcal{X})$ it must hold that $A[x+\delta]=b$ which implies $A \delta = 0$.
Consequently, all feasible $\delta\in \mathcal{D}(x; \mathcal{X})$ are in the null space of $A$.

Let $B = \begin{bmatrix}
    v_1 & \ldots & v_{|E|-r}
\end{bmatrix}$ form an orthonormal basis for the null space of matrix $A$, where $r$ is the rank of matrix $A$.
In the specific case for finding $\delta(x)$, c.f.~\eqref{eq:sfd}, $A$ has $N$ linearly independent rows and thus $r=N$.
To find $B$, we perform the singular value decomposition
\[
A = U \Sigma V^{\top},
\]
where $U\in \mathbb{R}^{N \times N}$ is a real orthogonal matrix, $\Sigma$ a diagonal matrix with singular values of $A$ on its diagonal in descending order and $V^{\top} \in \mathbb{R}^{|E| \times |E|}$ also is a real orthogonal matrix.
Then, $B$ corresponds to the last $|E| - r$ columns of $V$.

Having obtained orthonormal basis $B$, we can compute directional derivatives $\partial_{v_i} J(x)$ for each basis vector $v_i$, $i = 1, \dots, |E|-r$ in order to find \blue{the steepest feasible descent direction like the gradient but within the subspace spanned by the null space of $A$:}
\begin{align} \label{eq:subgrad}
    \delta'(x) = -\sum_{i=1}^{|E| - r} \partial_{v_i} J(x) \cdot v_i = -B \begin{bmatrix}
         \partial_{v_1} J(x) \\
         \vdots \\
         \partial_{v_{|E|-r}} J(x)
     \end{bmatrix},
\end{align}
for all $x \in \mathcal{X}_\varepsilon$ and normalizing the direction to unit length yields $\delta(x) = \delta'(x) / \|\delta'(x) \|_2$.

\subsection{SPSA via steepest feasible descent on $ \mathcal{X}^{\mathrm{eq}}$} \label{sec:spsa1}

Let us consider the following proxy of the steepest feasible descent direction for a function $J(x)$ in $x \in \mathcal{X}^{\mathrm{eq}}$:
\begin{align}
      \Tilde{\delta}(x) &=  \Tilde{\delta}(x, \Delta, \eta) \label{eq:res}\\
     &= \frac{J(x - \eta B \Delta) - J(x + \eta B \Delta)}{2\eta} B \Delta, \nonumber
\end{align}
where $B = \begin{bmatrix} v_1 & \dots & v_{|E|-r} \end{bmatrix}$, such that $\mathbb{C}\text{ol}(B) = \mathcal{N}(A)$, where $A$ as in \eqref{eq:A},
and perturbation vector $\Delta \in \mathbb{R}^{|E|-r}$.
Note that in the unconstrained case, $A$ is the null matrix and, consequently, $B = I$ so that we retrieve the expression in \eqref{eq:spsa_est}.

Consider the following assumptions:
\begin{enumerate}
    \item[\textbf{(A0)}] $J(x)$ is three times continuously differentiable with respect to $ x $ on $ \mathcal{X}^{\mathrm{eq}}$;
    \item[\textbf{(A1)}] $\Delta \sim \{-1,1\}^{|E|-r}$.
\end{enumerate}
The following theorem shows that, if these latter assumptions hold, then $\Tilde{\delta}(x)$ is an estimate of the steepest feasible descent direction $\delta'(x)$. 
\begin{theorem} \label{thm:sfd_est}
    { \em Let $\Tilde{\delta}(x)$ as in \eqref{eq:res}.
    Moreover, let us satisfy \textbf{(A0)} and \textbf{(A1)}. 
    Then, it holds that
    \begin{align}
        \mathbb{E}_{\Delta} \left[ \Tilde{\delta}(x)\right]  = \delta'(x) + \mathcal{O}(\eta^2).
    \end{align}
    }
    \begin{proof}
    See Appendix~\ref{sec:proof_sfd_est}.
    \end{proof}
\end{theorem}
Note that the estimation in Theorem~\ref{thm:sfd_est} simply ignores the norm constraint $\|\delta(x)\|_2 \leq 1$ in \eqref{eq:sfd} by considering $\delta'(x)$ instead of $\delta(x)$.

\subsection{Keeping SPSA on $ \mathcal{X}^{\mathrm{pos}}$} \label{sec:spsa2}
In the previous subsection, we have shown how the equality constraint can be incorporated via a base transformation and a steepest feasible descent direction staying on $ \mathcal{X}^{\mathrm{eq}}$ can be found.
In the following, we illustrate how to adjust \orange{the SPSA algorithm in \eqref{eq:spsa_rec}, so that $ x + \eta \Delta \in \mathcal{X}^{\mathrm{pos}}$:}
\begin{align} \label{eq:spsa}
    x(k+1) = \text{Proj}_{{\mathcal{X}^{(\mathrm{pos},\varepsilon)}}} \left[ x(k) - \alpha(k) G_k\left(x(k)\right)\right],
\end{align}
for $k\in \mathbb{N}$, where  $\text{Proj}_{Z}(x)$ is the projection of $x$ to point closest in the closed convex set $ Z $ in the $L^2$-norm \blue{ (viz., the element-wise projection $\text{Proj}_{{\mathcal{X}^{(\mathrm{pos},\varepsilon)}}}(x) = \max(\varepsilon, x)$)},
and perturbation scalar $\{\eta(k)\}_{k \in \mathbb{N}}$ in $G_k\left(x(k)\right)$ sufficiently small.
\blue{The latter can be ensured by setting $\eta(k) \leq \varepsilon $, for all $k \in \mathbb{N}$, as this ensures feasibility of perturbations in $G_k\left(x(k)\right)$, for any $\Delta \sim \{-1,1\}$; see \eqref{eq:spsa_est}.}
Effectively, the algorithm is applied on an interior set $\mathcal{X}^{(\mathrm{pos},\varepsilon)} \subset \mathcal{X}^{\mathrm{pos}}$, for some $\varepsilon$ small.
This approach can be extended for more generic constraints using numerical projections based on sequential quadratic programming \cite{shi2021sqp}.

As shown by \cite{Sadegh_1997}, one could apply a constrained SPSA on the original set ${\mathcal{X}}^{\mathrm{pos}}$.
To that end, introduce $\{\varepsilon(k)\} \to 0$ for $k \to \infty$, so that 
$
{\mathcal{X}}^{(\mathrm{pos}, \varepsilon(k))} \to {\mathcal{X}}^{\mathrm{pos}}  , 
\text{ for } k \to \infty.
$
Then, \orange{the algorithm in }\eqref{eq:spsa_rec} can be adapted as follows
\begin{align*}
    x(k+1) = \text{Proj}_{{\mathcal{X}^{\mathrm{pos}}}} \left[ x(k) - \alpha(k) G_k\left( \text{Proj}_{ \mathcal{X}^{(\mathrm{pos}, \varepsilon(k))}  } \left[ x(k)\right]\right)
    \right],
\end{align*}
for $k\in \mathbb{N}$, and is shown to converge a.s. \cite{Sadegh_1997}.

\subsection{The overall SPSA algorithm} \label{sec:spsa3}

\red{By integrating the approaches from \orange{s}ections~\ref{sec:spsa1} and~\ref{sec:spsa2}, we obtain the following SPSA algorithm:}
\begin{align} \label{eq:Alg}
    x(k+1) = \text{Proj}_{\mathcal{X_\varepsilon}} \left[x(k) + \alpha(k) \Tilde{\delta}_k(x(k))\right],
\end{align}
for $k \in \mathbb{N}$, and $\Tilde{\delta}_k(x(k)) = \Tilde{\delta}(x(k), \Delta(k), \eta(k))$, for $J(x) = S(P(x),C)$.
We simply assign $\varepsilon > 0$ an arbitrarily small value.
Furthermore, note that the operator $\text{Proj}_{\mathcal{X_\varepsilon}}[\cdot] $ can be easily implemented using the projection on a (scaled) simplex for each $x_{(i,\cdot)}$, for all $i \in V$, as in \cite{perez2020filtered}.

To state the convergence result for \eqref{eq:Alg}, the following additional assumptions for the sequences $\{\Delta (k)\}_{k \in \mathbb{N}}$, $\{\alpha (k)\}_{k \in \mathbb{N}}$, and $\{\eta (k)\}_{k \in \mathbb{N}}$ are sufficient:
\begin{enumerate}
    \item[\textbf{(A2)}] $\Delta(k) \sim \{-1,1\}^{|E|-r}$ is uniformly distributed for all $k$, $\{\alpha (k)\} =\{ \alpha / (\alpha_0 + k + 1)^{\gamma_{a}}\} $, and $\{\eta(k)\} = \{ \eta / (k + 1)^{\gamma_{\eta}}\} $, with $\alpha > 0$, $\alpha_0 \geq 0$, $\frac{1}{2} < \gamma_\alpha \leq 1$, and $\gamma_{\eta} > (1-\gamma_\alpha)/2$;
    \item[\textbf{(A3)}] $0 <\eta < \varepsilon / \sqrt{|E|-r}$, for $\varepsilon >0$.
\end{enumerate}
Note that \textbf{(A2)} extends the assumption on $\{\Delta (k)\}_{k \in \mathbb{N}}$ in \textbf{(A1)} to all iterations $k \in \mathbb{N}$.
Furthermore, \textbf{(A3)} ensures that for all $k \in \mathbb{N}$ and each realization $\Delta$, $x(k) \pm \eta B \Delta > 0$.
Consequently, $P(x(k) \pm \eta B \Delta)$ is an irreducible Markov chain and therefore $S(P(x(k)), C)$ is properly defined for all $k \in \mathbb{N}$.

\begin{theorem} \label{thm:est_sfd_conv}
{\em If assumptions \textbf{(A0)}, \textbf{(A2)} and \textbf{(A3)} are satisfied, the algorithm in \eqref{eq:Alg} converges almost surely \blue{to a stationary point of $S(P(x),C)$.}}

\begin{proof}
See Appendix~\ref{sec:proof_spsa}.
\end{proof}
\end{theorem}

\blue{An alternative approach for solving \eqref{eq:main_problem_non_rev} over $\mathcal{X}_\varepsilon$ would be to parameterize $x(u) = A^\dagger b + Bu \geq \varepsilon$, $u \in \mathbb {R}^ {|E| - r} $, where $A^\dagger = A^ {\top} (A A^\top)^ {-1}$ is the pseudo-inverse of $A$ as in \eqref{eq:A}, $B$ is an orthonormal basis for the null space of $A$, and $b$ is an appropriately sized vector of ones. 
Then, one can solve it using a constrained SPSA explained in Section~\ref{sec:spsa2}:
\begin{align}\label{eq:main_problem_non_rev_u}
    \begin{array}{rll}
    \min_{u \in \mathbb{R}^{|E|-r} } \hspace{1pt} &  S ( P(x(u)) , C ) \\
    \textup{s.t.} \hspace{1pt} & x(u) \geq \varepsilon.
    \end{array}
\end{align}
Note that applying SPSA on this problem formulation provides an equivalent steepest feasible descent estimator as in Theorem~\ref{thm:sfd_est}; however, to the best of our knowledge, there does not exist an off-the-shelf method for projecting on the constraint $x(u) = A^\dagger b + Bu \geq \varepsilon$.
One potential approach is to numerically solve the projection on the polytope $\{u \in \mathbb{R}^{|E|-r} : Bu \geq \varepsilon -A^\dagger b \}$, but this tends to be more computationally expensive compared to the well-studied projection on the probability simplex for projecting on $\mathcal{X_\varepsilon}$.}

\subsection{Constraining the stationary distribution} \label{sec:spsa4}

Theorem~\ref{thm:sfd_est} and Theorem~\ref{thm:est_sfd_conv} straightforwardly extend to the stationary distribution constrained case.
To that end, let
us find for feasible $x \in \mathcal{X}_\varepsilon^{\hat{\pi}} = \mathcal{X}_\varepsilon \cap \mathcal{X}^{\hat{\pi}}$,
\begin{align} \label{eq:sfd_pi}
    \delta(x) = \argmin_{\delta \in \mathcal{D} \left(x; \mathcal{X}^{\hat{\pi}}\right)} \ \partial_\delta S(P(x),C).
\end{align}
Then, let \blue{
\begin{align} \label{eq:A_stat}
    \begin{bmatrix}
    A^{\hat{\pi}} & b^{\hat{\pi}}
\end{bmatrix}  = \text{rref}\left(\begin{bmatrix}
A_1 & \bar{1}\\
A_2 & \Hat{\pi}
\end{bmatrix}\right),
\end{align}}
where the operator $\text{rref}\left( T \right)$ takes the row-reduced echelon form of a matrix $T$, $A_1$ is equal to the matrix in \eqref{eq:A}
and $A_2 = \begin{bmatrix}
    (A_2)_1 & \dots & (A_2)_N
\end{bmatrix}$,
such that $$(A_2)_i = \Hat{\pi}_i \cdot \begin{bmatrix} e_j \end{bmatrix}_{j : (i,j) \in E_i}, \quad \forall \, i \in V.$$
Using the row-reduced echelon form ensures that we omit any linear dependencies induced by constraining the stationary distribution.
It follows that \blue{$A^{\hat{\pi}} \in \mathbb{R}^{r^{\hat{\pi}} \times |E|}$ has $r^{\hat{\pi}}$ linearly independent rows and $b^{\hat{\pi}} \in \mathbb{R}^{r^{\hat{\pi}}}$.
Subsequently, we compute the orthonormal basis $B^{\hat{\pi}}$ for the null space of $A^{\hat{\pi}}$} and $\delta(x)$ using \eqref{eq:subgrad} and normalizing after.
Finally, the SPSA algorithm becomes, for $k \in \mathbb{N}$,
\begin{align} \label{eq:Alg_stat}
    x(k+1) = \text{Proj}_{\mathcal{X}_\varepsilon^{\hat{\hat{\pi}}}} \left[x(k) + \alpha(k) \Tilde{\delta}_k(x(k))\right]. 
\end{align}

Unlike in the previous setting, there is no straightforward closed-form projection on $\text{Proj}_{\mathcal{X_\varepsilon^{\hat{\pi}}}} [\cdot]$. 
However, we can use Dykstra's projection method of iteratively projecting on convex sets, see \cite{boyle1986method}, provided these projections can be carried out quickly.
To find $\text{Proj}_{\mathcal{X}_\varepsilon^{\hat{\pi}}}[x]$, let $\mathcal{X}_\varepsilon^{\hat{\pi}} = \mathcal{X}^1 \cap \mathcal{X}^2 \neq \emptyset$, with \blue{
\begin{align*} 
     &\mathcal{X}^{1} = \Big\{  x \in \mathbb{R}^{|E|} : \sum_{\ell \in E_i} x_{\ell} = 1, \ \forall \, i \in V , \hat{\pi} P(x) = \Hat{\pi} \Big\}, \\
    &\mathcal{X}^{2} = \left\{  x \in \mathbb{R}^{|E|} : \varepsilon \leq x \leq 1-\varepsilon \right\}.
\end{align*}}
One can project on $\mathcal{X}^{1}$ \blue{using some linear algebra, see Appendix~\ref{sec:app_proj}.}
The projection on $\mathcal{X}^{2}$ is straightforwardly done through:
\begin{align*}
    (\text{Proj}_{\mathcal{X}^{2}}[ x ])_{\ell} = \begin{cases}
        1-\varepsilon & \text{if } x_{\ell} > 1-\varepsilon, \\
        \varepsilon & \text{if } x_{\ell} < \varepsilon, \\
        x_{\ell} & \text{otherwise.} 
    \end{cases}
\end{align*}
\blue{Then, it is known that Dykstra's projection method (strongly) converges to $\text{Proj}_{\mathcal{X}^1 \cap \mathcal{X}^2} [x] = \text{Proj}_{X_\varepsilon^{\hat{\pi}}} [x]$, as proved in \cite{boyle1986method}.}

Note that strong convergence does not imply convergence in finitely many iterations (see \cite{bauschke2020dykstra} for a well-illustrated example).
For Dykstra's projection method, the number of iterations required to compute $\text{Proj}_{\mathcal{X}^1 \cap \mathcal{X}^2} [x]$ typically increases with the distance between $x$ and $\text{Proj}_{\mathcal{X}^1 \cap \mathcal{X}^2} [x]$. 
To ensure convergence within a reasonable number of iterations, one can bound the step size $\alpha(k) \Tilde{\delta}(x(k))$ in \eqref{eq:Alg} for each step $k$.
In all our numerical experiments, we saw a sufficiently rapid convergence, mainly due to choosing a relatively small step size $\alpha(k)$, see \orange{s}ections \ref{sec:ne} and \ref{sec:ne_2}.

\section{An extension of Markov chain connectivity to random graphs} \label{sec:ran}

In certain scenarios, graph structures are limited by the partial nature of our information. 
For instance, in transportation networks, the accessibility and availability of routes can fluctuate unpredictably due to extreme weather events \cite{papilloud2021vulnerability}, while power grids and communication networks exhibit analogous uncertainty, where the presence and reliability of the connection lines are subject to variability \cite{panteli2015modeling, soltan2015analysis, kubat1989estimation}.
In surveillance settings, the environment can be characterized by its susceptibility to changes, for example, through the potential obstruction of connections between various locations within the monitored area, often resulting from the presence of physical obstacles and impediments \cite{seder2005integrated,selek_Seder_Petrovic_2023}.

\subsection{Redistribution function}
We now introduce a framework to analyze the connectivity of a Markov chain on a graph whose vertex set is fixed, but the edge set is random.
\blue{We assume that only a subset of risky edges $\Tilde{E} \subset E$ are affected by uncertainty. 
Any risky edge $(i,j) \in \Tilde{E}$ fails (and therefore becomes inaccessible) with probability $q_{ij}\in [0,1]$.
Furthermore, we assume that $\red{\mathscr{G}}(V,E \setminus \Tilde{E})$ is strongly connected.
If a risky edge $(i,j)$ fails, it means that the original probability mass $P_{ij}$ must be redistributed over the remaining elements in the $i$th row of $P$.
Given a realization of accessible edges $\mathcal{E} \subseteq \Tilde{E}$, we calculate the adjusted Markov chain $Q(P, \mathcal{E}) = (Q_{ij}(P, \mathcal{E}))_{(i,j) \in V^2}$, where
\begin{align} \label{eq:adj_mc}
    Q_{ij}(P,\mathcal{E}) = \begin{cases}
    \frac{P_{ij}}{1 - \sum_{k:(i,k) \in \Tilde{E} \setminus \mathcal{E}} P_{ik}} & \text{ if } (i,j) \in \mathcal{E}; \\
    0 & \text{otherwise},
\end{cases}
\end{align}
for which it holds that $Q(P, \mathcal{E}) = P$ if $\mathcal{E}=\Tilde{E}$.
The \textit{expected Markov chain connectivity} of $P$ is then defined as
\begin{align} \label{eq:main_problem_stoch_derand}
 \mathbb{E}_{\mathcal{E}}\left[S ( Q(P, \mathcal{E}) , C )\right] = \sum_{i = 1}^{2^{|\Tilde{E}|}} \text{Prob}( \mathcal{E}_i ) \cdot S ( Q(P, \mathcal{E}_i) , C ),
\end{align}
where $\mathcal{E}_i \in 2^{\Tilde{E}}$ is a realization of the risky edges that are accessible.}  
In the specific case where edge failures are independent, the probability of observing a specific edge set $\mathcal{E}$ can be explicitly calculated as 
\begin{align} \label{eq:ex_prob}
    \text{Prob}( \mathcal{E} ) = \prod_{(i,j) \in \tilde{E} \setminus \mathcal{E}} q_{ij} \cdot \prod_{(i,j) \in \mathcal{E}} (1-q_{ij}).
\end{align}

\blue{To find a Markov chain $P$ with good connectivity properties while being robust in the uncertainty of the graph, we solve
\begin{align}\label{eq:main_problem_stoch}
\min_{P \in {\cal P} }\quad &  \mathbb{E}_{\mathcal{E}}\left[S ( Q(P, \mathcal{E}) , C )\right].
\end{align}
Note that enforcing the constraint $\pi(P) = \hat{\pi}$, i.e., 
\begin{align}\label{eq:main_problem_stoch_stat}
\min_{P \in {\cal P}^{\hat{\pi}} }\quad &  \mathbb{E}_{\mathcal{E}}\left[S ( Q(P, \mathcal{E}) , C )\right],
\end{align}
does not guarantee that the following constraint holds
\begin{align}
    \mathbb{E}_\mathcal{E} \left[ \pi(Q(P,\mathcal{E})) \right] = \hat{\pi}.
\end{align}}
\blue{As stated in the following theorem, we can obtain convexity for problem \eqref{eq:main_problem_stoch} (and problem \eqref{eq:main_problem_stoch_stat}) if we assume reversibility of $Q(P,\mathcal{E})$ together with the concurrent failure of reciprocal edges, that is, if edge $(i,j) \in \Tilde{E}$ fails, then $(j,i) \in \Tilde{E}$ also fails.}

\begin{theorem} \label{thm:convexity_rs}
    \blue{ \em  Let $C \in \mathbb{R}_{\geq 0}^{N \times N}$ be symmetric. Consider the stochastic optimization problem
    \begin{equation} \label{eq:main_problem_stoch_rev}
    \begin{aligned}
        \min_{x \in \mathcal{X}^{(\mathrm{sym}, 1)} }\quad &  \mathbb{E}_{\mathcal{E}}\left[S ( Q(x, \mathcal{E}) , C )\right], 
    \end{aligned}
    \end{equation}
    where $Q(x, \mathcal{E}) = Q(P(x), \mathcal{E})$.
    It follows that the reversibility condition \textbf{(R)} is satisfied for $P(x)$ via \eqref{eq:x_to_P} for all strongly connected $(\red{\mathscr{G}},x)$
    Moreover, let reciprocal edges fail concurrently.
    Then, \eqref{eq:main_problem_stoch_rev} is a convex problem.
    Additionally, if $C$ is a matrix with identical entries $C_{ij} = c >0$, \eqref{eq:main_problem_stoch_rev} is strictly convex.}
    
    \begin{proof} \blue{
    Assuming that the edges fail concurrently, we construct a new Markov chain where the edge weight of a failed edge is removed, see \eqref{eq:adj_mc}.
    Then, the optimization problem in \eqref{eq:main_problem_stoch_rev} reduces to 
    \begin{equation} \label{eq:prob_form_w}
    \begin{aligned}
        \min_{x \in \mathcal{X}^{(\mathrm{sym}, 1)}} \quad & \displaystyle \sum_{l = 1}^{2^{|\mathcal{E}|}} \text{Prob}( \mathcal{E}_l ) \cdot S ( Q(x, \mathcal{E}_l) , C ). \\
    \end{aligned}
    \end{equation}
    First, the feasible region $\mathcal{X}^{(\mathrm{sym}, 1)}$ is convex, being described only by linear constraints.
    We now show that $S ( Q(x, \mathcal{E}) , C )$ is a convex function in $x$, for all realizations $\mathcal{E}$.
    Given the symmetric failure of the edges, it follows that $Q(x, \mathcal{E})$ remains 
    reversible for any $\mathcal{E}$.
    Therefore, $S ( Q(x, \mathcal{E}) , C )$ is convex in $x$.
    
    It remains to show that if $C$ is a matrix with identical entries $C_{ij} = c >0$, \eqref{eq:prob_form_w} is strictly convex.
    This holds because, in this case, $S ( Q(x, \mathcal{E}) , C )$ is strictly convex in $x$ for all $\mathcal{E}$ and the sum of strictly convex functions is strictly convex.}
    \end{proof}
\end{theorem}

\blue{Using a similar argument as in Remark~\ref{rem:infty}, strongly connectedness in \eqref{eq:prob_form_w} is enforced via the objective function.}

\blue{Considering the limitations of reversible solutions demonstrated for the deterministic setting in Section~\ref{sec:por}, we will focus on exploring non-reversible solutions for \eqref{eq:main_problem_stoch} through
\begin{align}\label{eq:main_problem_stoch_non_rev}
    \min_{x \in {\cal X} } \quad &  \mathbb{E}_{\mathcal{E}}\left[S ( Q(x, \mathcal{E}) , C )\right],
\end{align}
where $Q(x, \mathcal{E}) = Q(P(x), \mathcal{E})$ with $P(x)$ via \eqref{eq:x_to_P2}, and for \eqref{eq:main_problem_stoch_stat}:
\begin{align}\label{eq:main_problem2_stoch_non_rev}
    \min_{x \in \mathcal{X}^{\hat{\pi} }} \quad &  \mathbb{E}_{\mathcal{E}}\left[S ( Q(x, \mathcal{E}) , C )\right].
\end{align}}
\red{Similarly to Section~\ref{sec:dif}, we solve these problems over $\mathcal{X}_\varepsilon$ and $\mathcal{X}_\varepsilon^{\hat{\pi}}$, respectively, which ensures that the SPSA algorithm remains applicable under the given constraints.}

\subsection{SPSA for Markov chains on random graphs} \label{sec:spsafmcorg}

When solving the problem formulated in \eqref{eq:main_problem_stoch}, an ideal scenario involves having access to perfect information regarding the probabilities of edge failure. 
However, this is not the case in many practical applications where edge failure probabilities must be inferred from the data.
In these cases, obtaining an exact analytical expression for $\mathbb{E}_{\mathcal{E}}\left[S ( Q(P, \mathcal{E}) , C )\right]$ as in \eqref{eq:main_problem_stoch_derand} becomes a difficult task.

Furthermore, the complexity of real-world systems often gives rise to intricate interdependencies between these failure probabilities. 
For instance, geographically close transmission power lines may have correlated failures due to extreme weather events, thereby introducing an additional layer of complexity into the modeling process \cite{neumayer2010network}. 
Accurate estimation of these correlations based on data is a notoriously difficult task, making the evaluation of $\mathbb{E}_{\mathcal{E}}\left[S ( Q(P, \mathcal{E}) , C )\right]$ problematic.

\red{We now show how SPSA is easily adapted to handle the random support setting, that is, we solve the stochastic optimization problem \eqref{eq:main_problem_stoch_non_rev}, 
where it is assumed that the distribution of $\mathcal{E}$ is unknown.}

Assume that the realizations of $\bm{\mathcal{E}} = (\mathcal{E}_1,\dots,\mathcal{E}_L)$ are subsequently presented to the decision maker.
In the following, we derive an estimator for the steepest feasible descent direction 
\begin{align*} 
    \delta_\mathcal{E}(x) = \argmin_{\delta \in \mathcal{D} \left(x; \mathcal{X}\right)} \ \partial_\delta \mathbb{E}_\mathcal{E}[ S(Q(x, \mathcal{E}), C) ],
\end{align*}
using graph realizations $\bm{\mathcal{E}}$.
To that end, first, observe that
\begin{align} \label{eq:exp_lim}
    \mathbb{E}_\mathcal{E}[ S(Q(x, \mathcal{E}), C) ] = \lim_{L\to \infty} \frac{1}{L} \sum_{l=1}^L S(Q(x, \mathcal{E}_l), C),
\end{align}
and for any feasible direction $\delta \in \mathcal{D}(x, \mathcal{X})$,
\begin{align*}
    \partial_\delta \mathbb{E}_\mathcal{E}[ S(Q(x, \mathcal{E}), C) ] = \lim_{L\to \infty} \frac{1}{L}  \sum_{l=1}^L \partial_\delta S(Q(x, \mathcal{E}_l), C),
\end{align*}
which follows from the independence between observations $\bm{\mathcal{E}}$, chosen $\delta$, and the finiteness of the collection of possible graph realizations.
It then follows that
\begin{align}
    \delta_\mathcal{E}'(x) &= -\sum_{i=1}^{|E|-r} \partial_{v_i} \left[ \lim_{L\to \infty}\frac{1}{L} \sum_{l=1}^L S(Q(x, \mathcal{E}_l), C)\right] \cdot v_i \nonumber\\
    &= -\sum_{i=1}^{|E|-r} \left[ \lim_{L\to \infty}\frac{1}{L} \sum_{l=1}^L \partial_{v_i} S(Q(x, \mathcal{E}_l), C) \right] \cdot v_i \nonumber \\
    &= -\sum_{i=1}^{|E|-r} \mathbb{E}_\mathcal{E}\left[  \partial_{v_i} S(Q(x, \mathcal{E}), C)\right] \cdot v_i, \label{eq:sfd_ran_subgrad}
\end{align}
where $v_i$ are the orthonormal basis vectors in $B=\begin{bmatrix} v_1 & \dots & v_{|E|-r} \end{bmatrix}$ obtained using the system of equations $Ax=b$ as in Section~\ref{sec:dif}.
\blue{Note that the constraint set $Ax=b$ is not affected by the realization $\mathcal{E}$ and, consequently, $B$ remains the same for all realizations of $\mathcal{E}$.}
Finally, $\delta_\mathcal{E}(x) = \delta_\mathcal{E}'(x) /\sum_{\ell \in E} (\delta_\mathcal{E}'(x))_\ell$.

Given a finite number of graph realizations $\bm{\mathcal{E}}$, define
\begin{align} \label{eq:j_L}
    \Bar{J}(x) = \Bar{J}(x, \bm{\mathcal{E}}) = \frac{1}{L}\sum_{l=1}^L S ( Q(x,\mathcal{E}_l) , C ),
\end{align}
which is an unbiased estimator for $\mathbb{E}_\mathcal{E}\left[  S(Q(x, \mathcal{E}), C)\right]$ for any $L \in \mathbb{N}$.
Now, let $B = \begin{bmatrix} v_1 & \dots & v_{|E|-r} \end{bmatrix}$, such that $\mathbb{C}\text{ol}(B) = \mathcal{N}(A)$, with $A$ as in \eqref{eq:A} and 
\begin{align}
 \Tilde{\delta}_{\mathcal{E}}(x)   = &  \ \Tilde{\delta}_{\mathcal{E}}(x, \Delta, \eta)  \nonumber \\
 = & \ \frac{\Bar{J}(x - \eta B \Delta) - \Bar{J}(x + \eta B \Delta)}{2 \eta} B\Delta, \label{eq:spsa_rs} 
\end{align}
such that $\Tilde{\delta}_\mathcal{E}(x)$ is an estimate of the steepest feasible descent direction $\delta_\mathcal{E}'(x)$ in \eqref{eq:sfd_ran_subgrad}, see Theorem~\ref{thm:est_sfd_ran}.

\begin{theorem} \label{thm:est_sfd_ran}
    { \em Let $\Tilde{\delta}_{\mathcal{E}}(x(k))$, for some $k$, as in \eqref{eq:spsa_rs}, where $\Bar{J}_L(x; \mathcal{E})$ as in \eqref{eq:j_L}, and assume \textbf{(A1)} and \textbf{(A3)}.
    Then, for any $L \in \mathbb{N}$, it holds that
    \begin{align}
        \mathbb{E}_{\Delta} \left[ \mathbb{E}_\mathcal{E} \left[ \Tilde{\delta}_{\mathcal{E}}(x)\right]\right]  = \delta_\mathcal{E}'(x) + \mathcal{O}(\eta^2).
    \end{align}
    }
    \begin{proof}
    The proof is analogous to that of Theorem~\ref{thm:sfd_est} using the fact that $\Bar{J}_L(x, \mathcal{E})$ is an unbiased estimator for $\mathbb{E}_\mathcal{E}\left[  S(Q(x, \mathcal{E}), C)\right]$. 
    \end{proof}
\end{theorem}

Theorem~\ref{thm:est_sfd_ran} allows for so-called \textit{online optimization}, \blue{where graph realizations are presented on a streaming basis.
In this context, one can make an iteration $k\in \mathbb{N}$, where $\Tilde{\delta}_{\mathcal{E}, k}(x(k)) = \Tilde{\delta}_{\mathcal{E}}(x(k), \Delta(k), \eta(k))$ based on an observed collection of realized edge sets $\bm{\mathcal{E}} = (\mathcal{E}_1, \ldots, \mathcal{E}_L)$.}
More formally, our policy $x(k)$ is updated via
\begin{align} \label{eq:Alg_ran}
    x(k+1) = \text{Proj}_{\mathcal{X_\varepsilon}} \left[x(k) + \alpha(k) \Tilde{\delta}_{\mathcal{E}, k}(x(k))\right], \,  \forall \, k\in \mathbb{N}.
\end{align}
The next theorem proves the convergence of this algorithm.

\begin{theorem} \label{thm:est_sfd_conv_ran}
{\em If assumptions 
\textbf{(A2)} and \textbf{(A3)} are satisfied, then the algorithm in \eqref{eq:Alg_ran} converges almost surely \blue{to a stationary point of $\mathbb{E}_{\mathcal{E}}\left[S ( Q(x, \mathcal{E}) , C )\right]$.}}

\begin{proof}
See Appendix~\ref{sec:proof_spsa_ran}.
\end{proof}
\end{theorem}

Note that enforcing the stationary distribution as in \eqref{eq:main_problem_stoch_stat} can be done similarly using the methods discussed in Section~\ref{sec:spsa4}.


\section{Numerical experiments}\label{sec:ne}
In this section, we demonstrate how we can numerically minimize $ S ( P(x) , C ) $ (or $\mathbb{E}_\mathcal{E}[ S(Q(x, \mathcal{E}), C) ]$ in the random support setting) on the feasible set $\mathcal{X}_\varepsilon^{\hat{\pi}}$.
For all experiments, we choose $\varepsilon = 10^{-4}$ and a single graph sample for our gradient estimator, that is $L=1$. Table~\ref{tab:param} provides an overview of the SPSA parameter settings. All numerical experiments were performed on a laptop with a 1.8 GHz processor and 12 GB RAM.

\begin{table}[!h]

\caption{SPSA parameters for the experiments of the next sections.}
\vspace{-0.25cm}
\label{tab:param}
\begin{center}
\begin{tabular}{|c||c|c|c|c|c|}
\hline
Exp. sec. & $\alpha$ & $\alpha_0$ & $\eta$ & $\gamma_\alpha$ & $\gamma_\eta$ \\ \hline
\ref{sec:scal} & $0.01$ & $100.000$ & $10^{-8}$ & $0.602$ & $0.200$ \\
\ref{sec:cor_fail} & $0.01$ & $100.000$ & $10^{-8}$ & $0.602$ & $0.200$ \\
\ref{sec:surv_fix} & $1$ & $500.000 $ & $10^{-8}$ & $0.602$ &  $0.200$ \\
\ref{sec:surv_ran} & $0.001$ & $50.000$  & $10^{-8}$ & $0.602$ &  $0.200$  \\ \hline
\end{tabular}
\end{center}

\end{table}

\subsection{Scalability of SPSA vs. IPOPT} \label{sec:scal}

In this subsection, we show the numerical results for SPSA for a network instance with $N=68$ nodes in which we progressively increase the size of the random support $\mathcal{E}$ and compare them with the Interior Point OPTimizer (IPOPT) solver as a benchmark.
IPOPT is a non-commercial solver capable of handling nonlinear (possibly non-convex) optimization problems with constraints ensuring primal feasibility while optimizing. 
It is readily available in the \texttt{cyipopt} package \cite{cyipopt}.

\begin{figure}
\vspace{-0.3cm}
    \centerline{
    \includegraphics[width=0.85\columnwidth]{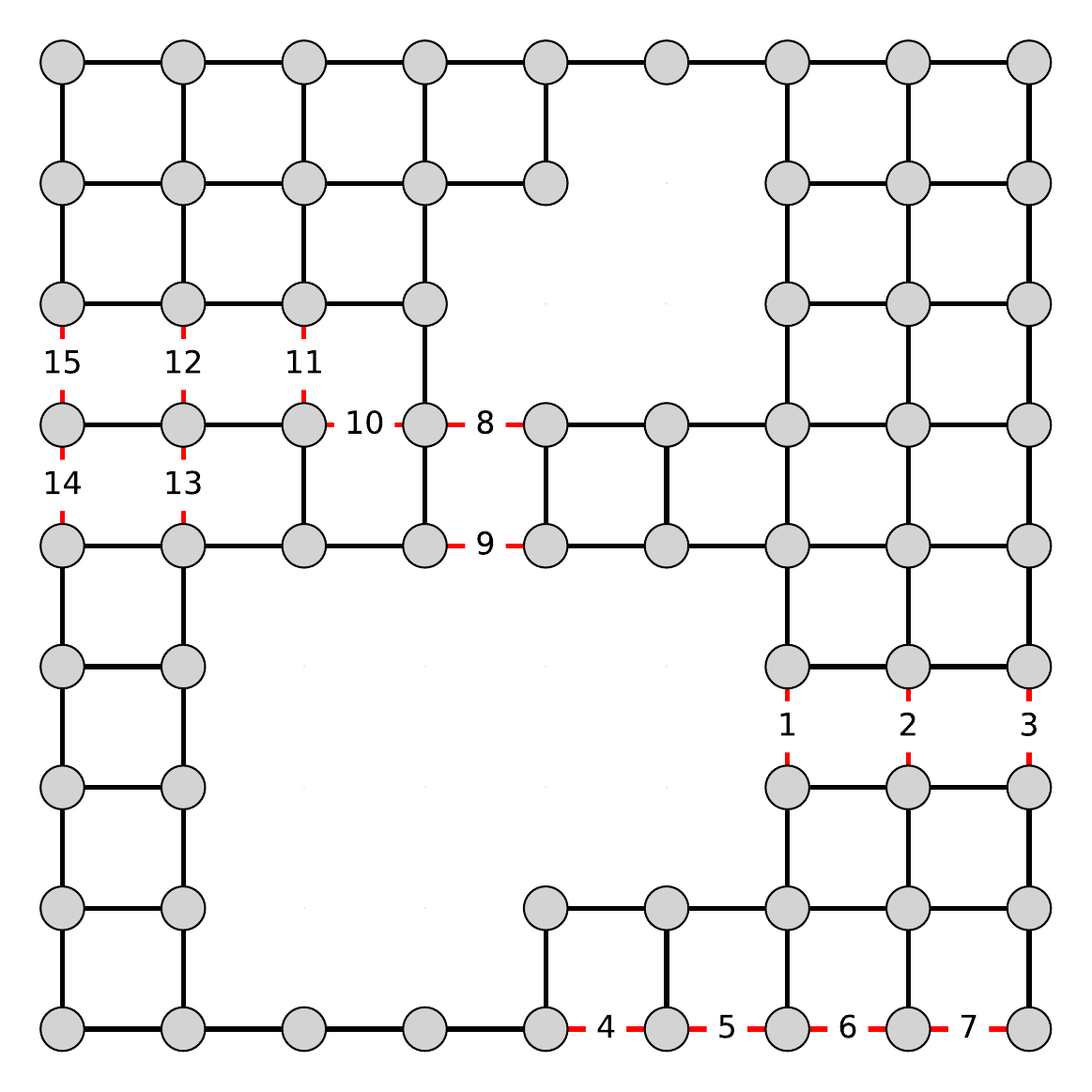}}
    \caption{Graph with $N = 68$ nodes similar to \cite{patel2015robotic}. 
    Red edges with labels indicate risky edges with a nonzero failure probability.}
    \label{fig:grid_n68}
\end{figure}

Consider the grid with bidirectional edges in Fig.~\ref{fig:grid_n68}. The edges in $\Tilde{E}$ that can fail (with probability $0.1$) are displayed in red and labeled.
It can be easily verified that $\red{\mathscr{G}} = (V,E \setminus \Tilde{E})$ is strongly connected.
Our goal is to find a reversible policy $P(x)$ minimizing $\mathbb{E}_\mathcal{E}\left[S ( Q(x, \mathcal{E}) , C^{\hat{\pi}} )\right]$, such that its stationary distribution is uniform, i.e., $\pi(P(x)) = \frac{1}{N}$.
Note that condition \textbf{(R)} implies that $P(x) = P(x)^{\top}$ and so we obtain the following problem:
\begin{align}\label{eq:problem_rev_unif}
\min_{x \in {\cal X}^{\mathrm{sym}} \cap {\cal X}_{\varepsilon}}\quad &  \mathbb{E}_\mathcal{E}\left[S ( Q(x, \mathcal{E}) , C^{\hat{\pi}} )\right],
\end{align}
which is strictly convex by Theorem~\ref{thm:convexity_rs}.
Note that in the current setting, we can parameterize the problem using a single variable $x_{ij}$ for each pair of bidirectional edges $(i,j), (j,i) \in E$.
Accordingly, we can write the equality constraints in ${\cal X}^{\mathrm{sym}} \cap {\cal X}_{\varepsilon}$ using $Ax=b$, where
\[
A_{il} = \begin{cases}
    1 & \text{if} \ i \ \text{is sender or receiver in edge $l$};\\
    0 & \text{otherwise},
\end{cases}
\]
such that $l$ corresponds to the $l$th edge in the lexicographic ordering of $E$, and $b = \bar{1}$.
We then apply the algorithm in \eqref{eq:Alg_ran} with $\text{Proj}_{{\cal X}^{\mathrm{sym}} \cap {\cal X}_{\varepsilon}}[x]$ utilizing Dykstra's algorithm (see Section~\ref{sec:spsa4}) and orthonormal basis $B$ in $\Tilde{\delta}_{\mathcal{E}}(x)$, such that $\mathbb{C}\text{ol}(B) = \mathcal{N}(A)$.

\blue{To compare our SPSA approach, we formulate an equivalent problem formulation of minimizing $\mathbb{E}_\mathcal{E}\left[S(Q(x(u), \mathcal{E}),C^{\hat{\pi}})\right]$, akin to \eqref{eq:main_problem_non_rev_u}.}
Then, we solve the problem using the IPOPT solver with a tolerance of $10^{-1}$. To start the IPOPT, we first collect a sample of 10 initial starting points and choose the one with the lowest objective value.
We provide the same initial starting point to both SPSA and IPOPT methods and terminate the SPSA algorithm if its objective value improves the IPOPT optimal solution, which we denote by $x^{\text{IPOPT}}$.
We evaluate the objective every $l = 50.000$ iterations by means of Polyak-Rupert averaging \cite{polyak1990new, ruppert1988efficient}, that is, we consider an average of the previous 50\% of observations;
\begin{align*}
    \bar{x}^{(l)}_i = \frac{1}{l\cdot i - \floor*{l\cdot i /2}} \sum_{k=\floor*{l\cdot i /2}}^{l\cdot i} x(k),
\end{align*}
and we terminate the algorithm if that solution improves the IPOPT solution, i.e., $| J(\bar{x}_i^{(l)}) - J(x^{\text{IPOPT}})| < 10^{-3}$, for $J(x) = \mathbb{E}_\mathcal{E}\left[S ( Q(x, \mathcal{E}) , C^{\hat{\pi}} )\right]$, or 
when a large number $K = 5 \cdot 10^6$ of iterations is achieved.
However, in all our experiments, convergence was reached in less than $K$ steps.

\begin{figure}
\centerline{
    \includegraphics[width=0.925\columnwidth]{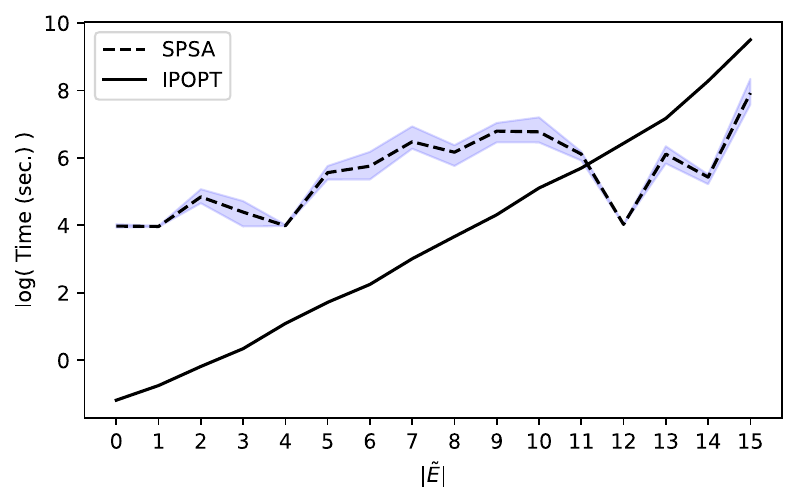}}
    \caption{Running time (in log(sec)) for SPSA and IPOPT for increasing size $\Tilde{E}$ with range of running time depicted in blue.}
    \label{fig:cal_E_inc}
\end{figure}

We calibrate SPSA according to \eqref{eq:Alg} with the settings as in Table~\ref{tab:param} and repeat the simulation for $|\Tilde{E}| = 0,1,\ldots,15$, where $\Tilde{E}$ is the set of edges with the lowest labels included in Fig.~\ref{fig:grid_n68}.
We show the running time for IPOPT and the SPSA mean and range for 10 repetitions for each $|\Tilde{E}| = 0,1,\ldots,15$ in Fig.~\ref{fig:cal_E_inc}.
\blue{Note that the situation $|\Tilde{E}| = 0$ corresponds to the fixed support case in which IPOPT is faster.
However, our analysis reveals that the running time of IPOPT demonstrates a pattern closely resembling exponential growth. }
In contrast, SPSA is not subject to a similar exponential increase in computational demand. 
Specifically, SPSA tends to underperform in smaller instances, which may be attributed to inadequate tuning of hyperparameters.
Nevertheless, it should be noted that the running time of SPSA consistently remains below 1.2 hours across all tested instances.
A comparison of SPSA with other solvers, such as the Sequential Least Squares Quadratic Programming (SLSQP) solver from the \texttt{Scipy} package, gave similar results to those presented in Fig.~\ref{fig:cal_E_inc}.

\subsection{Correlated edge failures} \label{sec:cor_fail}

This subsection considers the case of correlated failures, which means that $\mathbb{E}_\mathcal{E}\left[S ( Q(x,\mathcal{E}) , C )\right]$ cannot be easily computed assuming the independence of edge failure.
Therefore, we cannot solve problem \eqref{eq:problem_rev_unif} efficiently using objective evaluations as the probabilities for appearing graphs cannot be calculated assuming independent edge failure as in \eqref{eq:ex_prob}.
Moreover, using a sample average approximation typically causes scalability problems, especially when the set of failing edges $\Tilde{E}$ becomes large.
These cases are particularly favorable for SPSA, which has no issues dealing with edge failure correlations, as we can simply observe a sample of the correlated distribution of graphs at each iteration for gradient estimation.

In a second experiment, we solve the instance in which the set of failing edges $\Tilde{E}$ includes the edges with labels 1, 2, 3, 8, and 9 in Fig.~\ref{fig:grid_n68}. 
We assume that each edge fails with probability $q= 0.5$ and the joint failure of edges is strongly correlated (0.85).
Let $\Bar{J}(x)$ denote a sample-average approximation of $\mathbb{E}_\mathcal{E}\left[S ( Q(x,\mathcal{E}) , C )\right]$ using a large sample of $10.000$ graphs following~\cite{chen2015generating}.
We ran the same SPSA algorithm as in Section~\ref{sec:scal} with parameters as in Table~\ref{tab:param} and terminated the algorithm when $| \Bar{J}(\bar{x}_i^{(l)}) - \Bar{J}(\bar{x}^{(l)}_{i-1})| < 10^{-3}$, which we verify every $l=50.000$ iterations.
The result is $\Bar{J}(x^{\text{cor}}) = 344,08$.

We compare the correlated solution with the uncorrelated solution found using IPOPT (with tolerance $10^{-4}$), that is, we assume in this setting i.i.d.~edge failures with probability 0.5, $ \Bar{J}(x^{\text{uncor}}) = 350,17$.
Note that the SPSA algorithm is able to find a solution that is lower than the (estimated) objective by roughly 1.74\%.

\section{Application: robotic surveillance} \label{sec:ne_2}

In this section, we validate the applicability of our method in the context of surveillance, in which problem (\ref{eq:main_problem2}) has a natural interpretation.
Here, $\hat{\pi}_i$ is the anticipated probability that an intruder appears at a node $ i$ (e.g., uniform or based on historical data), so that $ P(x) $ with $ \pi ( P(x) ) = \hat{\pi} $ that minimizes $ S (P(x), C^{\hat{\pi}} ) $ (or $\mathbb{E}_\mathcal{E}\left[S(Q(x, \mathcal{E}), C^{\hat{\pi}})\right]$ in the random support setting), is a surveillance policy with a minimal average number of steps to capture the intruder. 
Specifically, this surveillance problem with $ C =  C^{\hat{\pi}} $ has been studied by \cite{patel2015robotic}, in which the authors find the optimal surveillance policy with minimal Kemeny constant assuming reversibility of $P(x)$ in the fixed support setting.

In the following, we first consider a surveillance problem similar to that in~\cite{patel2015robotic} and then an analogous problem in the random support setting.

\subsubsection{Fixed support} \label{sec:surv_fix}
\begin{figure}
    
    \centerline{
    \includegraphics[width=0.85\columnwidth]{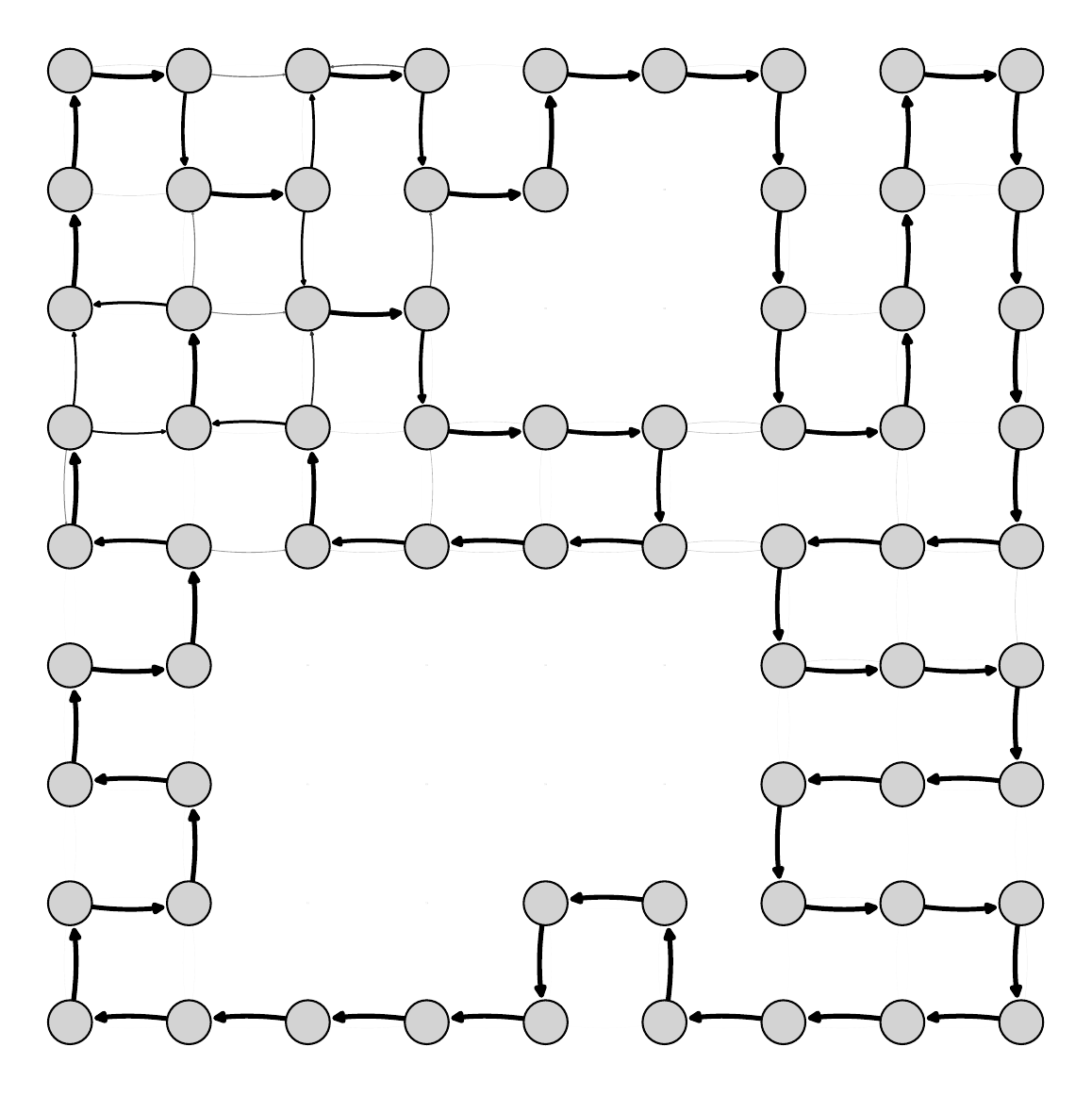}}
    \caption{Solution to \eqref{eq:prob1} for the instance in Fig.~\ref{fig:grid_n68}. The thickness of each arrow is proportional to the corresponding transition probability.}
    \label{fig:grid_n68_sol}
  \end{figure}

Consider again the graph with $N=68$ nodes shown in Fig.~\ref{fig:grid_n68}. 
Similarly to \cite{patel2015robotic}, we formulate an optimization problem to find a stochastic surveillance policy for an agent that visits each node with equal probability and minimizes the Kemeny constant, namely
\begin{equation} \label{eq:prob1}
    \begin{array}{rll}
    \displaystyle \min_{x \in \mathcal{X}_\varepsilon^{\hat{\pi}}} & S\left(P(x), C^{\hat{\pi}}\right),
    \end{array}
\end{equation}
where $\hat{\pi} = \frac{1}{N}$. 
We apply SPSA as detailed in Section~\ref{sec:spsa4}.
We terminate the algorithm if the condition $$| S(P(x_i^{(l)}), C^{\hat{\pi}} )- S(P(x_{i-1}^{(l)}), C^{\hat{\pi}}) )| < 10^{-3}$$ is met, which we verify every $l = 1.000$ iterations.
For the parameters as in Table~\ref{tab:param}, the algorithm terminated in 254 seconds and we present the resulting Markov chain in Fig.~\ref{fig:grid_n68_sol}.
Note that the solution is only locally optimal and approximates a Hamiltonian cycle.

To test the validity of a solution in a surveillance context, we assume that 500 intruders appear uniformly at random in a network node and reside there for 45 time units. For simplicity, we assume that individual intrusions appear sequentially at time $0, 45, 90, \dots$.
We initialize the surveillance agent uniformly at random on the graph, and after each time unit, the surveillance agent instantaneously transitions to an adjacent node according to the chosen surveillance policy $P(x)$. 
If the surveillance agents and the intruder are in the same node, the intruder is caught.

We compare our non-reversible solution $P(x^{\text{non-rev}})$ with 
the Markov chain with minimal Kemeny constant assuming reversibility $P(x^{\text{rev}})$ \cite{patel2015robotic}.
We report the results for 500 simulations in Table \ref{tab:fixed_surv}, for each of the \orange{two} solutions. 
As in the experiment by \cite{patel2015robotic}, we see that the number of intruders caught increases with a decrease in the Kemeny constant of the surveillance policy; the policy with the lowest Kemeny constant catches the most intruders.
Inspecting the $P(x^{\text{non-rev}})$ policy in Fig.~\ref{fig:grid_n68_sol}, we see that a strong level of (clockwise) directionality is achieved (note that the Hamiltonian cycle is the extreme case of directionality), which the $P(x^{\text{rev}})$ policy cannot obtain.
This directionality feature is crucial for our approach to obtain lower MFPTs and, thus, higher intruder detection rates. 

\begin{table}[!h]
\caption{Intruders caught (in \%) and average connectivity value achieved by optimal reversible and non-reversible stochastic surveillance policies.}
\vspace{-0.25cm}
\label{tab:fixed_surv}
\begin{center}
\begin{tabular}{|c||c|c|c|c||c|}
\hline
\multicolumn{6}{|c|}{\rule{0pt}{2ex}Fixed support} \\ \hline %
Sol. & Min (\%) & Mean (\%) & Max (\%) & SD & Obj. value \\ \hline
rev & 19,00 & 26,36 & 32,40 & 1,94    & 192,7           \\
\textbf{non-rev}   & \textbf{50,40} & \textbf{57,31} & \textbf{64,80} & \textbf{2,33}    & \textbf{51,8}  \\ \hline
\multicolumn{6}{|c|}{\rule{0pt}{2ex}Random support} \\ \hline
Sol. & Min (\%) & Mean (\%) & Max (\%) & SD & Obj. value \\ \hline
rev & 17,80 & 23,72 & 31,00 & 2,02 & 266,3 \\
\textbf{non-rev}   & \textbf{23,00} & \textbf{33,70} & \textbf{41,20} & \textbf{3,30}    & \textbf{162,4}  \\ \hline
\end{tabular}
\end{center}
\end{table}

\subsubsection{Random support} \label{sec:surv_ran}

\blue{To show that non-reversible solutions outperform reversible ones in the random support case, }we now consider an instance of the graph depicted in Fig.~\ref{fig:grid_n68} in which a subset of \blue{risky} edges can be inaccessible with some probability.

To find a stochastic surveillance policy, we solve the following problem
\begin{equation} \label{eq:prob_temp}
    \begin{array}{rll}
    \displaystyle \min_{x \in \mathcal{X}_\varepsilon^{\hat{\pi}}} & \mathbb{E}_\mathcal{E} \left[S\left(Q(x, \mathcal{E}), C^{\Hat{\pi}}\right) \right],
    \end{array}
\end{equation}
where we aim to visit every node with equal probability in case there is no failure in the graph ($\hat{\pi} = \frac{1}{N}$).
The problem objective is closely related to the expected Kemeny constant; however, note that while $C^{\Hat{\pi}}$ is fixed, the true expected Kemeny constant is $C^{\pi(Q(P, \mathcal{E}))}$, which depends on the graph realization.

In this example, we compare the globally optimal reversible policy $P(x^{\text{rev}})$, where $\Tilde{E}$ is again the set of failing edges corresponding to labels 1, 2, 3, 8, and 9, with failure probability $0.5$, 
with the locally optimal non-reversible policy $P(x^{\text{non-rev}})$.
To find $P(x^{\text{rev}})$, we use the IPOPT algorithm with tolerance $10^{-4}$.
For finding $P(x^{\text{non-rev}})$, we apply SPSA in \eqref{eq:Alg_ran}, 
with $\text{Proj}_{\mathcal{X_\varepsilon^{\hat{\pi}}}} \left[x\right]$, and orthonormal basis $B$ in $\Tilde{\delta}_{\mathcal{E}}(x)$ spanning the null space of $A$ as in \eqref{eq:A_stat}.
We terminate the algorithm if condition $$| \mathbb{E}_\mathcal{E}[S(Q(x_{i-1}^{(l)}, \mathcal{E}), C^{\hat{\pi}} )] - \mathbb{E}_\mathcal{E}[S(Q(x_{i-1}^{(l)}, \mathcal{E}), C^{\hat{\pi}}) )]| < 10^{-3}$$ is met, which we verify every $l=10.000$ iterations.
The remaining hyperparameters are reported in Table~\ref{tab:param}.

To test the validity of our solution in a surveillance context, we performed $500$ simulations with sampled graphs $\mathcal{G}_l = (V, \mathcal{E}_l)$, for $l=1,\ldots, 500$, with edge failure probability $0.5$ for every $(i,j) \in \mathcal{E}$.
As in the previous experiment, 500 intruders appear uniformly at random in a network node and reside there for 45 time units. 
We initialize the surveillance agent uniformly at random on the graph at time $t=0$ and after each time unit, the surveillance agent instantaneously transitions to an adjacent node according to the surveillance policy $Q(x,\mathcal{E}_l)$ 
We report the results for these 500 simulations in Table \ref{tab:fixed_surv}, for the globally optimal reversible solution $P(x^{\text{rev}})$ and locally optimal non-reversible solution $P(x^{\text{non-rev}})$. 
Analogously to the previous experiment with fixed support, we see that the number of intruders caught increases with a lower objective value for the surveillance policy; the non-reversible policy outperforms the reversible policy. 

\section{Conclusion and future work}
\label{sec:conclusion}
We have shown how to find Markov chains minimizing a weighted sum of mean first passage times, which generalizes existing metrics in the literature used in the connectivity analysis of the corresponding network. 
Previous studies strongly relied on reversible Markov chains (or undirected graphs), aiming to formulate and solve a convex problem. 
Our work shows that a strictly better solution for minimizing the sum of mean first passage times can be found when the reversibility constraint is dropped, hence considering Markov chain optimization on directed graphs.
To solve this more general class of problems, we extended SPSA to Markov chain optimization that does not require the reversibility constraint.
Furthermore, we showed how SPSA can efficiently solve the minimization of weighted mean first passage times in the case of random and potentially correlated failure of network edges.
Finally, we evaluated a locally optimal non-reversible Markov chain as a surveillance policy and showed that it outperforms reversible Markov chain policies, both when edges are safe, as well as when edge failure occurs.

The generalized SPSA algorithm that we propose can be applied to other problems that involve optimization on the probability simplex, such as portfolio optimization \cite{jondeau2006optimal}, and other resource allocation problems \cite{patriksson2008survey}.
In the specific context of surveillance, future work will explore an online version of the surveillance policy optimization in which the empirical frequency of intruder presence is progressively learned rather than being assumed to be a uniform stationary distribution.

\appendix

\subsection{Infeasibility SPSA} \label{sec:inf_spsa}
To illustrate why standard SPSA does not work, let us consider the following Markov chain on the complete graph $\red{\mathscr{G}}$ with $N=3$ nodes and stochastic matrix
\[
P(x) = \begin{bmatrix} 
        0 & \frac{1}{2} & \frac{1}{2} \\[2pt]
        \frac{1}{2} & 0 & \frac{1}{2}\\[2pt]
        \frac{1}{2} & \frac{1}{2} & 0
    \end{bmatrix}.
\]
A standard application of SPSA for $J(x) = S(P(x), C)$, for some $C \in \mathbb{R}_{\geq 0}^{N \times N}$ and realization $\Delta = (1,1,-1,1,1,1)$ yields
\[
P(x + \eta \Delta) = \begin{bmatrix} 
        0 & \frac{1}{2} + \eta & \frac{1}{2} +\eta\\[2pt]
        \frac{1}{2} - \eta & 0 & \frac{1}{2} + \eta\\[2pt]
        \frac{1}{2}+ \eta & \frac{1}{2}+ \eta & 0
    \end{bmatrix} \not \in \mathcal{P},
\]
for any small $\eta >0$. It follows that $\Pi$ diverges and therefore $D$, $M$, and $S(P(x + \eta \Delta),C)$ (and $S(P(x - \eta \Delta),C)$) is not well defined.

\subsection{Projection on $\mathcal{X}_1$} \label{sec:app_proj}
\blue{ We illustrate here how to project on 
\[
\mathcal{X}^{1} = \Big\{  x \in \mathbb{R}^{|E|} : \sum_{\ell \in E_i} x_{\ell} = 1, \ \forall \, i \in V , \hat{\pi} P(x) = \Hat{\pi} \Big\},
\]
which spans a linear subspace in $\mathbb{R}^{|E|}$. We can perform an orthogonal projection on the subspace $\mathcal{X}_1$ for a given point $x$ via
\begin{align*}
\begin{split}
    &\text{Proj}_{\mathcal{X}^{1}}[ x ] = \\
    & \, = B^{\hat{\pi}} \left((B^{\hat{\pi}})^{\top} B^{\hat{\pi}}\right)^{-1} (B^{\hat{\pi}})^{\top} \left(x - (A^{\hat{\pi}})^{\dagger} b^{\hat{\pi}}\right) + (A^{\hat{\pi}})^{\dagger} b^{\hat{\pi}} \\
    & \, = B^{\hat{\pi}} (B^{\hat{\pi}})^{\top} \left(x - (A^{\hat{\pi}})^{\dagger} b^{\hat{\pi}}\right) + (A^{\hat{\pi}})^{\dagger} b^{\hat{\pi}},
\end{split}
\end{align*}
where $(A^{\hat{\pi}})^\dagger = (A^{\hat{\pi}})^{\top} (A^{\hat{\pi}} (A^{\hat{\pi}})^\top)^{-1}$ is the pseudo-inverse of $A^{\hat{\pi}}$ and the latter equality follows from the fact that $(B^{\hat{\pi}})^{\top} B^{\hat{\pi}} = \left((B^{\hat{\pi}})^{\top} B^{\hat{\pi}}\right)^{-1} = I$, as $B^{\hat{\pi}}$ contains an orthonormal basis.}

\subsection{Proof of Lemma~\ref{lem:surj}} \label{sec:proof_surj}

\blue{Observe that if the weighted graph $(\red{\mathscr{G}},x)$, where $x\in \mathcal{X}^{\mathrm{sym}}$, is strongly connected, the Markov chain $P(x)$ via \eqref{eq:x_to_P} satisfies \textbf{(R)}, i.e., $P(x) \in \mathcal{P}^{\mathrm{rev}}$ \cite{ghosh2008minimizing}.
Moreover, it has stationary distribution \orange{$\pi$, where for all $i \in V$} \begin{align} \label{eq:stat_un_G}
    \pi_i(P(x)) = \frac{\sum_{j \in E_i} x_{(i,j)}}{\sum_{(i,j) \in E} x_{(i,j)}}.
\end{align}
Now, we show that any $P \in \mathcal{P}^{\mathrm{rev}}$ has a representation $x \in \mathcal{X}^{(\mathrm{sym}, 1)}$, such that $P = P(x)$ where $P(x)$ via \eqref{eq:x_to_P}.
Indeed, it is enough to take $x_{(i,j)} = \pi_i(P) P_{ i j }$ for all $(i,j) \in E$ so that $P(x) \in \mathcal{P}^{\mathrm{rev}}$ via \eqref{eq:x_to_P}, which shows that $P(x)$ is surjective.
More specifically, we can show this representation is unique.
To that end, use \eqref{eq:stat_un_G} and the constraint $\sum_{(i,j) \in E}x_{(i,j)} = 1$ to observe that any $x \in \mathcal{X}^{(\mathrm{sym}, 1)}$ has fixed values
$\pi_i(P(x)) = \sum_{j \in E_i} x_{(i,j)}$, for all $i \in V$, as follows from the fact that any candidate $x\in \mathcal{X}^{(\mathrm{sym}, 1)}$ must have $\pi(P(x)) = \pi(P)$.
As a result, $x_{(i,j)}$  must be uniquely determined, as shown by $$P_{ij}(x) =\frac{x_{(i,j)}}{\sum_{j \in E_i} x_{(i,j)}} = P_{ij}.$$ 
This uniqueness implies that $x$ is a unique representation for $P$ via $P(x)$. Moreover, since $x$ is strongly connected as follows from the irreducibility of $P \in \mathcal{P}^{\mathrm{rev}}$, the mapping $P(x)$ is a bijection from $\mathcal{X}^{(\mathrm{sym}, 1)} \cap \{x \in \mathcal{X}^{\mathrm{sym}} : (\red{\mathscr{G}},x) \text{ is strongly connected}\}$ to $\mathcal{P}^{\mathrm{rev}}$.} 

\subsection{Proof of Proposition~\ref{prop:tour}} \label{sec:proof_tour}

\red{The optimization problem~\eqref{eq:opt_P_nonrev} has been studied in detailed in \cite{Borkar_Miclo_2020}, in which the authors proved that the optimal solution $\Hat{P}$ is a Hamiltonian cycle, see Theorem 2 therein.} The same paper also shows that $\sum_{i \neq j} M_{ij} = \frac{N(N-1)}{2}$ for all $i$, from which it easily follows that
\[
    S ( \Hat{P},  C^{\bar{1}}  ) = \sum_{i=1}^N \sum_{j \neq i} M_{ij} = N \frac{N(N-1)}{2} = \frac{N^3 - N^2}{2}.
\]

\subsection{Proof of Proposition~\ref{prop:ul_bound}} \label{sec:proof_ub_lb}

Consider the family of complete undirected graphs $(\red{\mathscr{G}},x)$ with $x \in \mathcal{X}^{\mathrm{sym}}$.
The complete graph with $N$ nodes, with $x_{\ell} = 1$, for all $\ell \in E$, has been shown to have total effective graph resistance $R_{tot}(x) = N-1$, see \cite{ellens2011effective}.
Using this fact in combination with \eqref{eq:egr} and \eqref{eq:obj} provides the lower bound in \eqref{eq:lu_bound}, that is
\[
S(\hat{P}, C^{\bar{1}}) = N(N-1) (N-1) = N^3 - 2N^2 + N.
\]

The upper bound in \eqref{eq:lu_bound} follows from the fact that the Hamiltonian graph with maximal $R_{tot}(x)$ is a circular graph (as adding edges to $\red{\mathscr{G}}$ cannot increase $R_{tot}(x) $, \cite{ghosh2008minimizing}). In \cite{lukovits1999resistance}, it has been shown that for a circular graph $(\red{\mathscr{G}},x)$ of size $N$, $R_{tot}(x)  = \frac{N^3 - N}{12}$, where $x_{\ell} = 1$ if $\ell \in E$ is part of the circle and $x_{\ell} = 0$ otherwise.
Since $(\red{\mathscr{G}},x)$ is a circular graph, we have $\sum_{\ell \in E} x_{\ell} = N$.
Rewriting \eqref{eq:egr} yields:
\[
S(\Hat{P}, C^{\bar{1}}) = 2N \left(\frac{N^3 - N}{12}\right) = \frac{N^4 - N^2}{6}.
\]

\subsection{Proof of Theorem~\ref{thm:sfd_est}} \label{sec:proof_sfd_est}

Assume $x \in \mathcal{X_\varepsilon}$ and let $J(u)$ denote the linearly (re)parameterized version of $J(x(u))$, where $P(x(u)) = \sum_{\ell=(i,j) \in E} x_\ell(u) e_i e_j^{\top}$ if $x(u) \geq 0$, with $e_i$ denoting the $i$th (column) basis vector, and $x(u) = A^\dagger b + Bu$, such that $u \in\mathbb {R}^ {|E| - r} $ and $b$ as in \eqref{eq:A}.
Expanding $J( u \pm \eta \Delta)$, for $u = B^\dagger(x - A^\dagger b ) = B^{\top}(x - A^\dagger b )$, yields
\orange{\begin{align*}
    J(u \pm \eta \Delta)& = J(u) \pm \eta \sum_{i=1}^{|E|-r} \partial_{u_i} J(u) \Delta_i  \\
    & \hspace{11pt} + \frac{\eta^2}{2} \sum_{i,j=1}^{|E|-r} \frac{\partial J(u)}{\partial_{u_i} \partial_{u_j}} \Delta_i \Delta_j \\
    &  \hspace{11pt} \pm \frac{\eta^3}{6} \sum_{i,j,k=1}^{|E|-r} \frac{\partial J(u)}{\partial_{u_i} \partial_{u_j} \partial_{u_k}} \Delta_i \Delta_j \Delta_k + \mathcal{O}(\eta^4).
\end{align*}}
\red{We can use the expansion to estimate partial derivatives $\partial_{u_l} J(u)$, for any $l$, as we show in the following.} Basic algebra shows that \red{for a given $\Delta \in \{-1,1\}^d$, and some $ y, y' \in [ u - \eta \Delta , u + \eta \Delta ]$;}  
\blue{\red{
\begin{align*}
\begin{split}
    &\hspace{-15pt} \frac{J(u + \eta \Delta) - J(u - \eta \Delta)}{2 \eta \Delta_l} =  \\
    & = \sum_{i=1}^{|E|-r} \partial_{u_i} J(u) \frac{\Delta_i}{\Delta_l} \\
    & \hspace{11pt} + \frac{\eta^2}{6} \sum_{i,j,k=1}^{|E|-r} \hspace{-2pt}\big( \frac{\partial J(y)}{\partial_{u_i} \partial_{u_j} \partial_{u_k}}-  \frac{\partial J(y')}{\partial_{u_i} \partial_{u_j} \partial_{u_k}}\big)  \hspace{-1pt} \frac{\Delta_i \Delta_j \Delta_k}{\Delta_l}.
\end{split}
\end{align*}}}
Note that by \textbf{(A0)} \blue{it holds that} the 3rd order derivative is bounded (for a proof, use the fact that $ \mathcal{X}$ is compact together with Weierstrass theorem). 
\red{Then, taking the expectation over $\Delta \sim \text{Rad}^d$ (i.e., each $\Delta_i$ i.i.d.\ uniform over $\{-1,1\}$) by assumption \textbf{(A1)}, we obtain:
\begin{align*}
    \mathbb{E}_\Delta\left[\frac{J(u + \eta \Delta) - J(u - \eta \Delta)}{2 \eta \Delta_l}\right] = \partial_{u_l} J(u) + \mathcal{O}(\eta^2),
\end{align*}}
\red{i.e., effectively estimating $\partial_{u_l} J(u)$.} Now using that $x(u) = A^\dagger b + Bu$, and so
\begin{align*}
\begin{split}
    x( u \pm \eta \Delta) &= A^{\dagger}b + B[u \pm \eta \Delta] = x(u) \pm \eta B \Delta,
\end{split}
\end{align*}
\blue{and $\partial_{u_l} x(u) = \partial_{v_l} Bu = v_l$, we obtain
\begin{align*}
    \mathbb{E}_\Delta\left[\frac{J(x + \eta B \Delta) - J(x - \eta B \Delta)}{2 \eta \Delta_l}\right] = \partial_{v_l} J(x) + \mathcal{O}(\eta^2).
\end{align*}}
Combining all elements in vector form \blue{
\begin{align}
    & \sum_{l=1}^{|E|-r} -v_l \cdot \mathbb{E}_\Delta\left[\frac{J(x+ \eta B\Delta) - J(x-\eta B\Delta)}{2 \eta \Delta_l}\right]+ \mathcal{O}(\eta^2) \nonumber \\
    & =\sum_{l=1}^{|E|-r} v_l \cdot \mathbb{E}_\Delta\left[\frac{J(x -  \eta B\Delta) - J(x + \eta B \Delta)}{2 \eta \Delta_l}\right]+ \mathcal{O}(\eta^2) \nonumber \\
    &= \mathbb{E}_\Delta \left[ \frac{J(x - \eta B \Delta) - J(x + \eta B \Delta)}{2\eta} B \Delta \right] + \mathcal{O}(\eta^2) \nonumber \\
    &= \mathbb{E}_\Delta\left[\Tilde{\delta}(x)\right] + \mathcal{O}(\eta^2) = \delta'(x), \label{eq:exp_bias}
\end{align}}%
where we use in the second equality that $\Delta_i = \Delta^{-1}_i$, for all $i$ and $\Delta$, by \textbf{(A1)}.

\subsection{Proof of Theorem~\ref{thm:est_sfd_conv}} \label{sec:proof_spsa}

The proof follows a basic stochastic approximation argumentation assuming the projected ODE \blue{on a convex and bounded set $X$ given by}:
\begin{align*} 
\partial \text{Proj}_X(x_t) / \partial t = \lim_{h \to 0} \text{Proj}_{X}(x_t - h g(x)) / h,
\end{align*}
\blue{where $\{x_t : t \geq 0\}$ is the trajectory \red{of} the ODE, and $g(x)$ is the gradient of some function $J(x)$. Let $\Tilde{g}(x)$ be a (biased) gradient estimator of $g(x)$, so that
\begin{align*}
    \mathbb{E} \left[ \Tilde{g}(x(k)) \right] &= \mathbb{E} \left[ \Tilde{g}(x(k)) | x(0), \ldots, x(k-1) \right] \\ &= g(x(k)) + \beta(x(k)).
\end{align*}}
The following conditions from \cite{Kushner_Yin_2003}, see Section 5.2.1, can be applied to ensure convergence:
\begin{enumerate}
    \setlength\itemsep{1.5pt}
    \item[\textbf{(C1)}] \blue{$\sup_k \mathbb{E} [\|\Tilde{g}(x(k))\|_2^2] < \infty$;}
    \item[\textbf{(C2)}] $g(x)$ is continuous;
    \item[\textbf{(C3)}] $\sum_{k=0}^\infty \alpha(k) = \infty$, $\sum_{k=0}^\infty \alpha(k)^2 < \infty$;
    \item[\textbf{(C4)}] \blue{$\sum_{k=0}^\infty  \alpha(k) \| \beta(x(k)) \|_2 < \infty$,} w.p. 1;
    \item[\textbf{(C5)}] $J(x)$ and the constraints imposing $\mathcal{X}_\varepsilon$ are twice continuously differentiable.
\end{enumerate}

Conditions \textbf{(C1)} -- \textbf{(C5)} are also sufficient in the case where $g(x)$ is a projected gradient on a smooth manifold or, equivalently, the steepest descent direction.
In our setting, let $J(x) = S(P(x), C)$, and $X = \mathcal{X}_\varepsilon$.
Given that $\mathcal{X}_\varepsilon$ is smooth, we can consider $g(x) = -\delta'(x)$ as in \blue{\eqref{eq:subgrad}} to be the projected (negative)  gradient of some smooth surrogate function $H(x)$, for which $H(x) = J(x)$, for all $x \in \mathcal{X}_\varepsilon$, but it is defined on $\mathcal{X}_\varepsilon \cup \mathcal{X}_\varepsilon^{\perp}$\blue{, where $\mathcal{X}_\varepsilon^{\perp}$ is the orthogonal complement of $\mathcal{X}_\varepsilon$}.
In that setting, 
\begin{align*}
\delta'(x) = \lim_{h \to 0} \text{Proj}_{\mathcal{X}_\varepsilon}(x_t - h \nabla H(x)) / h.
\end{align*}

Now, assume that
\begin{align*}
    \mathbb{E}_\Delta \left[ \Tilde{\delta}(x(k)) \right] &= \mathbb{E}_\Delta \left[ \Tilde{\delta}(x(k)) | x(0), \ldots, x(k-1) \right] \\ &= \delta'(x(k)) + \beta(x(k)),
\end{align*}
where $\beta(x(k)) = \mathbb{E}[ \Tilde{\delta}(x(k)) ] - \delta'(x(k)) = \mathcal{O}(\eta^2)$ indicates the bias, see \eqref{eq:exp_bias}.
Observe that $J(x)$ is differentiable everywhere and therefore uniformly Lipschitz continuous with Lipschitz constant $\kappa$ by the smoothness of $J$ and compactness of $\mathcal{X}_\varepsilon$. 
\blue{To show that \textbf{(C1)} is satisfied, observe that $|J(x(k) \pm \eta(k) B\Delta(k)) - J(x(k))| < \infty$ as follows from choosing $\eta(k)$ via \textbf{(A3)}, so that $P(x(k) \pm \eta(k) B\Delta(k)) > 0$ and so is irreducible.
Then, $\Tilde{\delta}(x(k))$ can be bounded by}
\begin{align}
    |\Tilde{\delta}_i(x(k))| &\leq \frac{|J(x(k) + \eta(k) B\Delta(k)) - J(x(k))|}{2 \eta(k) |\Delta_i(k)|} \nonumber \\
    & \hspace{10pt} + \frac{|J(x(k) - \eta(k) B\Delta(k)) - J(x(k))|}{2 \eta(k) |\Delta_i(k)|} \nonumber \\
    & \leq \frac{\| \eta(k) B \Delta(k)\|_2 \kappa}{2 \eta(k) |\Delta_i(k)|} + \frac{\|- \eta(k) B \Delta(k)\|_2 \kappa}{2 \eta(k) |\Delta_i(k)|} \nonumber \\
    &= \frac{\|\Delta(k)\|_2}{|\Delta_i(k)|} \kappa =\kappa \sqrt{|E| - r}, \label{eq:lip_bound}
\end{align}
for $i = 1, \ldots, |E|$ and all $k\in \mathbb{N}$, where the latter equality follows from \blue{\textbf{(A2)}}. 
Condition \textbf{(C2)} follows from the fact that $J(x) = S(P(x),C)$ is smooth.
Condition \textbf{(C3)} is easily verified by basic calculus for the assumed sequence $\{\alpha (k)\} =\{ \alpha / (\alpha_0 + k + 1)^{\gamma_{a}}\} $, for $\frac{1}{2} < \gamma_{a} \leq 1$.
Similarly, \textbf{(C4)} is verified using the fact that the bias $\beta(x(k)) = \mathcal{O}(\eta^2)$ is fixed and \blue{$\sum_{k=0}^\infty  \alpha(k) \| \beta(x(k)) \|_2 < \infty$} for $\{\eta(k)\} = \{ \eta / (k + 1)^{\gamma_{\eta}}\} $, where $\gamma_{\eta} > (1-\gamma_a)/2$. 
Finally, \textbf{(C5)} follows from $J(x) = S(P(x),C)$ smooth and the affine constraints imposing $\mathcal{X}_\varepsilon$.

It remains to show for which $\eta$ it holds that $x \pm \eta B \Delta \in \mathcal{X}$, for all $\Delta \in \{-1,1\}^{|E|-r}$ and all $x \in \mathcal{X}_\varepsilon$. 
Given that $x \pm B \Delta$ preserves the equality constraints in $\mathcal{X}$, it is enough that $0 < \eta B \Delta < \varepsilon$ for all $\Delta \in \{-1,1\}^{|E|-r}$.
To this end, observe that $\|\Delta\|_2 = \sqrt{|E|-r}$ and so $\|B\Delta\|_2 =\sqrt{|E|-r}$, since $B$ is an orthonormal basis.
It follows that $\eta \|B \Delta \|_2 < \varepsilon$ implies $\eta < \varepsilon / \|B \Delta \|_2 = \eta < \varepsilon / \sqrt{|E| - r}$, see \textbf{(A3)}.

\subsection{Proof of Theorem~\ref{thm:est_sfd_conv_ran}} \label{sec:proof_spsa_ran}

We verify conditions \textbf{(C1)} -- \textbf{(C5)} as in the proof of Theorem~\ref{thm:est_sfd_conv}. Let $J(x) = \mathbb{E}_\mathcal{E}[S(Q(x, \mathcal{E}), C)]$, $g(x) = \delta_\mathcal{E}'(x)$ as in \eqref{eq:sfd_ran_subgrad}.
Similarly to the proof of Theorem~\ref{thm:est_sfd_conv}, assume that
\begin{align*}
    \mathbb{E}_\Delta \hspace{-1pt} \left[ \mathbb{E}_\mathcal{E} \hspace{-1pt} \left[ \Tilde{\delta}_{\mathcal{E}}(x(k)) \right]\right] \hspace{-1pt} &= \mathbb{E}_\Delta \hspace{-1pt} \left[ \mathbb{E}_\mathcal{E} \hspace{-1pt} \left[\Tilde{\delta}_{\mathcal{E}}(x(k) | x(0), \ldots, x(k-1) \right]\right] \\ &= \delta_\mathcal{E}'(x(k)) + \beta(x),
\end{align*}
where $\beta(x(k))$ is a random variable.
Condition \textbf{(C1)} is verified given that the graph collection is finite, $S(Q(x, \mathcal{E}), C)$ is smooth everywhere and therefore uniformly Lipschitz for all $\bm{\mathcal{E}} = (\mathcal{E}_1,\dots,\mathcal{E}_L)$, $L \in \mathbb{N}$.
This can be bounded similarly to that in \eqref{eq:lip_bound}.
Furthermore, \textbf{(C2)} easily follows from inspecting \eqref{eq:sfd_ran_subgrad}.
Conditions \textbf{(C3)} and \textbf{(C4)} follow the same arguments as shown in the proof of Theorem~\ref{thm:est_sfd_conv}, and \textbf{(C5)} follows from $\delta_\mathcal{E}'(x)$ smooth and the affine constraints imposing $\mathcal{X}_\varepsilon$.

\bibliographystyle{IEEEtran}
\bibliography{main_v4}

\begin{thebibliography}{10}
\providecommand{\url}[1]{#1}
\csname url@samestyle\endcsname
\providecommand{\newblock}{\relax}
\providecommand{\bibinfo}[2]{#2}
\providecommand{\BIBentrySTDinterwordspacing}{\spaceskip=0pt\relax}
\providecommand{\BIBentryALTinterwordstretchfactor}{4}
\providecommand{\BIBentryALTinterwordspacing}{\spaceskip=\fontdimen2\font plus
\BIBentryALTinterwordstretchfactor\fontdimen3\font minus \fontdimen4\font\relax}
\providecommand{\BIBforeignlanguage}[2]{{%
\expandafter\ifx\csname l@#1\endcsname\relax
\typeout{** WARNING: IEEEtran.bst: No hyphenation pattern has been}%
\typeout{** loaded for the language `#1'. Using the pattern for}%
\typeout{** the default language instead.}%
\else
\language=\csname l@#1\endcsname
\fi
#2}}
\providecommand{\BIBdecl}{\relax}
\BIBdecl

\bibitem{Hunter2016}
J.~J. Hunter, ``The computation of key properties of {M}arkov chains via perturbations,'' \emph{Linear Algebra and its Applications}, vol. 511, pp. 176--202, 2016.

\bibitem{Hunter2018}
------, ``The computation of the mean first passage times for {M}arkov chains,'' \emph{Linear Algebra and its Applications}, vol. 549, pp. 100--122, 2018.

\bibitem{Chou2014}
T.~Chou and M.~R.~D. Orsogna, ``First passage problems in biology,'' in \emph{First-Passage Phenomena and Their Applications}.\hskip 1em plus 0.5em minus 0.4em\relax {World} {Scientific}, 2014, pp. 306--345.

\bibitem{Kalantar2018}
N.~Kalantar and D.~Segal, ``Mean first-passage time and steady-state transfer rate in classical chains,'' \emph{The Journal of Physical Chemistry C}, vol. 123, no.~2, pp. 1021--1031, 2019.

\bibitem{Weiss2007}
G.~H. Weiss, ``First passage time problems in chemical physics,'' in \emph{Advances in Chemical Physics}.\hskip 1em plus 0.5em minus 0.4em\relax J.~Wiley {\&} Sons, 2007, pp. 1--18.

\bibitem{klein1993resistance}
D.~J. Klein and M.~Randi{\'c}, ``Resistance distance,'' \emph{Journal of {M}athematical {C}hemistry}, vol.~12, pp. 81--95, 1993.

\bibitem{ellens2011effective}
W.~Ellens, F.~M. Spieksma, P.~Van~Mieghem, A.~Jamakovic, and R.~E. Kooij, ``Effective graph resistance,'' \emph{Linear Algebra and its Applications}, vol. 435, no.~10, pp. 2491--2506, 2011.

\bibitem{bianchi2019kirchhoffian}
M.~Bianchi, J.~L. Palacios, A.~Torriero, and A.~L. Wirkierman, ``Kirchhoffian indices for weighted digraphs,'' \emph{Discrete Applied Mathematics}, vol. 255, pp. 142--154, 2019.

\bibitem{kemeny1976finite}
J.~G. Kemeny and J.~L. Snell, \emph{Finite {M}arkov {C}hains {W}ith a {N}ew {A}ppendix "{G}eneralization of a {F}undamental {M}atrix"}.\hskip 1em plus 0.5em minus 0.4em\relax Springer, 1976.

\bibitem{ghosh2008minimizing}
A.~Ghosh, S.~Boyd, and A.~Saberi, ``Minimizing effective resistance of a graph,'' \emph{SIAM Review}, vol.~50, no.~1, pp. 37--66, 2008.

\bibitem{patel2015robotic}
R.~Patel, P.~Agharkar, and F.~Bullo, ``Robotic surveillance and {M}arkov chains with minimal weighted {K}emeny constant,'' \emph{IEEE Transactions on Automatic Control}, vol.~60, no.~12, pp. 3156--3167, 2015.

\bibitem{Duan_Bullo_2020}
X.~Duan and F.~Bullo, ``Markov chain--based stochastic strategies for robotic surveillance,'' \emph{Annual Review of Control, Robotics, and Autonomous Systems}, vol.~4, no.~1, pp. 243--264, 2021.

\bibitem{spall1992multivariate}
J.~C. Spall, ``Multivariate stochastic approximation using a simultaneous perturbation gradient approximation,'' \emph{IEEE {T}ransactions on {A}utomatic {C}ontrol}, vol.~37, no.~3, pp. 332--341, 1992.

\bibitem{Sadegh_1997}
P.~Sadegh, ``Constrained optimization via stochastic approximation with a simultaneous perturbation gradient approximation,'' \emph{Automatica}, vol.~33, no.~5, p. 889–892, 1997.

\bibitem{berkhout2019analysis}
J.~Berkhout and B.~F. Heidergott, ``Analysis of {M}arkov influence graphs,'' \emph{Operations Research}, vol.~67, no.~3, pp. 892--904, 2019.

\bibitem{aldous2014reversible}
D.~Aldous and J.~A. Fill, ``Reversible {M}arkov {C}hains and {R}andom {W}alks on {G}raphs,'' 2002, unfinished monograph, recompiled 2014, available at \url{https://www.stat.berkeley.edu/users/aldous/RWG/book.pdf}.

\bibitem{Zocca2021}
A.~Zocca and B.~Zwart, ``Optimization of stochastic lossy transport networks and applications to power grids,'' \emph{Stochastic Systems}, vol.~11, no.~1, pp. 34--59, 2021.

\bibitem{wang2015network}
X.~Wang, Y.~Ko{\c{c}}, R.~E. Kooij, and P.~Van~Mieghem, ``A network approach for power grid robustness against cascading failures,'' in \emph{2015 7th {I}nternational {W}orkshop on {R}eliable {N}etworks {D}esign and {M}odeling (RNDM)}.\hskip 1em plus 0.5em minus 0.4em\relax IEEE, 2015, pp. 208--214.

\bibitem{dorfler2010synchronization}
F.~Dorfler and F.~Bullo, ``Synchronization of power networks: Network reduction and effective resistance,'' \emph{IFAC Proceedings Volumes}, vol.~43, no.~19, pp. 197--202, 2010.

\bibitem{tizghadam2010betweenness}
A.~Tizghadam and A.~Leon-Garcia, ``Betweenness centrality and resistance distance in communication networks,'' \emph{IEEE Network}, vol.~24, no.~6, pp. 10--16, 2010.

\bibitem{rueda2017robustness}
D.~F. Rueda, E.~Calle, and J.~L. Marzo, ``Robustness comparison of 15 real telecommunication networks: Structural and centrality measurements,'' \emph{Journal of Network and Systems Management}, vol.~25, pp. 269--289, 2017.

\bibitem{yang2018designing}
C.~Yang, J.~Mao, X.~Qian, and P.~Wei, ``Designing robust air transportation networks via minimizing total effective resistance,'' \emph{IEEE Transactions on Intelligent Transportation Systems}, vol.~20, no.~6, pp. 2353--2366, 2018.

\bibitem{yilmaz2020kemeny}
S.~Yilmaz, E.~Dudkina, M.~Bin, E.~Crisostomi, P.~Ferraro, R.~Murray-Smith, T.~Parisini, L.~Stone, and R.~Shorten, ``Kemeny-based testing for {COVID}-19,'' \emph{PLOS ONE}, vol.~15, no.~11, p. e0242401, 2020.

\bibitem{Borkar_Miclo_2020}
V.~Borkar and L.~Miclo, ``On the fastest finite {M}arkov processes,'' \emph{Journal of Mathematical Analysis and Applications}, vol. 481, no.~2, p. 123488, 2020.

\bibitem{shi2021sqp}
J.~Shi and J.~C. Spall, ``{SQP}-based projection {SPSA} algorithm for stochastic optimization with inequality constraints,'' in \emph{2021 American control conference (ACC)}.\hskip 1em plus 0.5em minus 0.4em\relax IEEE, 2021, pp. 1244--1249.

\bibitem{heidergott2003taylor}
B.~Heidergott and A.~Hordijk, ``Taylor series expansions for stationary {M}arkov chains,'' \emph{Advances in Applied Probability}, vol.~35, no.~4, pp. 1046--1070, 2003.

\bibitem{leder2010approximation}
N.~Leder, B.~Heidergott, and A.~Hordijk, ``An approximation approach for the deviation matrix of continuous-time {M}arkov processes with application to {M}arkov decision theory,'' \emph{Operations {R}esearch}, vol.~58, no. 4-part-1, pp. 918--932, 2010.

\bibitem{perez2020filtered}
G.~Perez, M.~Barlaud, L.~Fillatre, and J.-C. R{\'e}gin, ``A filtered bucket-clustering method for projection onto the simplex and the $\pmb {l}_1$ ball,'' \emph{Mathematical Programming}, vol. 182, no. 1-2, pp. 445--464, 2020.

\bibitem{boyle1986method}
J.~P. Boyle and R.~L. Dykstra, ``A method for finding projections onto the intersection of convex sets in {H}ilbert spaces,'' in \emph{Advances in Order Restricted Statistical Inference: Proceedings of the Symposium on Order Restricted Statistical Inference held in Iowa City, Iowa, September 11--13, 1985}.\hskip 1em plus 0.5em minus 0.4em\relax Springer, 1986, pp. 28--47.

\bibitem{bauschke2020dykstra}
H.~H. Bauschke, R.~S. Burachik, D.~B. Herman, and C.~Y. Kaya, ``On {D}ykstra’s algorithm: finite convergence, stalling, and the method of alternating projections,'' \emph{Optimization Letters}, vol.~14, pp. 1975--1987, 2020.

\bibitem{papilloud2021vulnerability}
T.~Papilloud and M.~Keiler, ``Vulnerability patterns of road network to extreme floods based on accessibility measures,'' \emph{Transportation research part D: transport and environment}, vol. 100, p. 103045, 2021.

\bibitem{panteli2015modeling}
M.~Panteli and P.~Mancarella, ``Modeling and evaluating the resilience of critical electrical power infrastructure to extreme weather events,'' \emph{IEEE Systems Journal}, vol.~11, no.~3, pp. 1733--1742, 2015.

\bibitem{soltan2015analysis}
S.~Soltan, D.~Mazauric, and G.~Zussman, ``Analysis of failures in power grids,'' \emph{IEEE Transactions on Control of Network Systems}, vol.~4, no.~2, pp. 288--300, 2015.

\bibitem{kubat1989estimation}
P.~Kubat, ``Estimation of reliability for communication/computer networks simulation/analytic approach,'' \emph{IEEE Transactions on Communications}, vol.~37, no.~9, pp. 927--933, 1989.

\bibitem{seder2005integrated}
M.~Seder, K.~Macek, and I.~Petrovic, ``An integrated approach to real-time mobile robot control in partially known indoor environments,'' in \emph{31st Annual Conference of IEEE Industrial Electronics Society, 2005. IECON 2005.}\hskip 1em plus 0.5em minus 0.4em\relax IEEE, 2005, pp. 6--pp.

\bibitem{selek_Seder_Petrovic_2023}
A.~Šelek, M.~Seder, and I.~Petrović, ``Smooth autonomous patrolling for a differential-drive mobile robot in dynamic environments,'' \emph{Sensors}, vol.~23, no.~17, p. 7421, 2023.

\bibitem{neumayer2010network}
S.~Neumayer and E.~Modiano, ``Network reliability with geographically correlated failures,'' in \emph{2010 Proceedings IEEE INFOCOM}.\hskip 1em plus 0.5em minus 0.4em\relax IEEE, 2010, pp. 1--9.

\bibitem{cyipopt}
J.~K. Moore \emph{et~al.}, ``{cyipopt},'' \url{https://cyipopt.readthedocs.io/en/stable/}, 2024, version 1.1.0.

\bibitem{polyak1990new}
B.~T. Polyak, ``New stochastic approximation type procedures,'' \emph{Avtomatika i Telemekhanika}, vol.~7, no. 98-107, p.~2, 1990.

\bibitem{ruppert1988efficient}
D.~Ruppert, ``Efficient estimations from a slowly convergent {R}obbins-{M}onro process,'' Cornell University Operations Research and Industrial Engineering, Tech. Rep., 1988.

\bibitem{chen2015generating}
M.~Chen, ``Generating nonnegatively correlated binary random variates,'' \emph{The Stata Journal: Promoting Communications on Statistics and Stata}, vol.~15, no.~1, p. 301–308, 2015.

\bibitem{jondeau2006optimal}
E.~Jondeau and M.~Rockinger, ``Optimal portfolio allocation under higher moments,'' \emph{European Financial Management}, vol.~12, no.~1, pp. 29--55, 2006.

\bibitem{patriksson2008survey}
M.~Patriksson, ``A survey on the continuous nonlinear resource allocation problem,'' \emph{European Journal of Operational Research}, vol. 185, no.~1, pp. 1--46, 2008.

\bibitem{lukovits1999resistance}
I.~Lukovits, S.~Nikoli{\'c}, and N.~Trinajsti{\'c}, ``Resistance distance in regular graphs,'' \emph{International Journal of Quantum Chemistry}, vol.~71, no.~3, pp. 217--225, 1999.

\bibitem{Kushner_Yin_2003}
H.~Kushner and G.~G. Yin, \emph{Stochastic Approximation and Recursive Algorithms and Applications}, 2nd~ed., ser. Stochastic Modelling and Applied Probability.\hskip 1em plus 0.5em minus 0.4em\relax Springer Science \& Business Media, 2003, vol.~35.

\end{thebibliography}

\end{document}